\documentclass[]{interact}
\usepackage[numbers,sort&compress]{natbib}
\bibpunct[, ]{[}{]}{,}{n}{,}{,}
\makeatletter
\def\NAT@def@citea{\def\@citea{\NAT@separator}}
\makeatother

\usepackage{amsmath,amsfonts,amssymb,latexsym}
\usepackage{amsthm}
\usepackage{mathtools}
\usepackage[normalem]{ulem}
\usepackage{enumerate,longtable,color,multirow,mathrsfs}
\usepackage[]{ccaption}
\captiondelim{ } 
\usepackage[bookmarks]{hyperref}

\usepackage[linesnumbered,ruled,vlined]{algorithm2e}

\numberwithin{equation}{section}

\theoremstyle{plain}
\newtheorem{theorem}{Theorem}
\newtheorem{problem}{Open problem}

\newtheorem{lemma}{Lemma}[section]
\newtheorem{corollary}[lemma]{Corollary}
\newtheorem{proposition}[lemma]{Proposition}
\newtheorem{proposition and definition}[lemma]{Proposition \& Definition}

\theoremstyle{definition}
\newtheorem{definition}[lemma]{Definition}
\newtheorem{example}[lemma]{Example}
\theoremstyle{remark}
\newtheorem{remark}[lemma]{Remark}
\newtheorem{notation}[lemma]{Notation}

\newcommand{\ad}{\mathrm{ad}}
\newcommand{\af}{\mathfrak{aff}}

\newcommand{\C}{\mathbb C}
\newcommand{\D}{\mathfrak{D}}

\newcommand{\co}{\coloneqq}
\newcommand{\diam}{\mathrm{diam}}
\newcommand{\dist}{\mathrm{dist}}
\newcommand{\Der}{\mathrm{Der}}
\newcommand{\id}{\mathrm{Id}}
\newcommand{\Ker}{\mathrm{Ker}}
\newcommand{\pr}{\operatorname{pr}}
\newcommand{\N}{\mathbb N}
\newcommand{\Zz}{\mathcal{Z}}
\newcommand{\R}{\mathbb R}
\newcommand{\F}{\mathbb F}
\newcommand{\G}{\mathcal G}

\newcommand{\s}{\mathrm{span}}
\newcommand{\Sc}{\mathcal S}
\newcommand{\Nn}{\mathcal{N}}

\newif\ifbibdoi \bibdoitrue
\providecommand{\bibdoi}[1]{\ifbibdoi\ doi:\href{https://doi.org/#1}{{\urlstyle{same}\nolinkurl{#1}}}.\fi}

\begin{document}

\title{Linear structure graphs associated with commuting graphs of Lie algebras}

\author{
	\name{Vu A. Le\textsuperscript{a,b}\thanks{Corresponding author: Vu A. Le. Email: \underline{vula@uel.edu.vn}. \\ \indent $\,$ Regarding the use of AI assistance, please refer to the declaration in Subsection~\ref{Subsec5.3New}.}, Hoa Q. Duong\textsuperscript{c}
	and Tuan A. Nguyen\textsuperscript{d}}
	\affil{\textsuperscript{a}University of Economics and Law, Ho Chi Minh City, Vietnam;
	\textsuperscript{b}Vietnam National University, Ho Chi Minh City, Vietnam;
	\textsuperscript{c}University of Finance - Marketing, Ho Chi Minh City, Vietnam;
	\textsuperscript{d}Ho Chi Minh City University of Education, Ho Chi Minh City, Vietnam
	}
}

\maketitle

\begin{abstract}
The commuting graph (CG) of a finite-dimensional Lie algebra over a field of characteristic zero forgets the linear structure, and a non-trivial center creates redundant vertices. We introduce linear structure graphs (LSGs), whose vertices are the nonzero vectors of a vector space and whose edge set is stable under linear operations. Each vertex generates a finite ascending sequence of subspaces, which defines its LSG-height; this height governs connectivity. We characterize LSGs of height at most 1 and show that the CG of a Lie algebra with trivial center has height at most 1 if and only if the algebra is a commutative transitive Lie algebra (CTLA). The reduced commuting graph (RCG) restricts the CG to a complement of the center; it is well defined up to LSG-isomorphism, preserves the connectivity, diameter and number of components of the CG, and is compatible with direct sums. Finally, we characterize almost abelian Lie algebras (AALAs) by their RCGs.
\end{abstract}

\begin{keywords}
	Linear structure graph (LSG); LSG-monomorphism; LSG-isomorphism; Lie algebra; LSG-adjacency; LSG-height; LSG-closed subspace; commuting graph (CG); reduced commuting graph (RCG); almost abelian Lie algebra (AALA)
\end{keywords}

\begin{amscode}
	05C25; 17B05; 05C12; 05C40; 05C69; 05C70; 05C76; 17B20; 17B30
\end{amscode}

\tableofcontents

\section*{Introduction}

Associating a graph with an algebraic structure and recovering algebraic information from the combinatorial invariants of that graph has become a standard approach in modern algebra. This approach goes back to the work of Brauer and Fowler on centralizers of involutions in finite groups~\cite{BrauerFowler1955}, which preceded the formal notion of commuting graphs and later played a role in the classification of finite simple groups, and to Beck's zero-divisor graphs of commutative rings~\cite{Beck1988}. Among such graphs, the \emph{commuting graphs} (CGs) occupy a central place: their vertices are the elements lying outside the center, and two distinct vertices are adjacent if and only if  they commute. CGs have been studied extensively for groups, rings and matrix algebras (see, e.g.,~\cite{Akbari2008, Dorbidi2024, Miguel2016, Nam2024}); their connectivity, diameter and clique structure reflect essential features of the underlying algebraic objects.

Much of this research is driven by the \emph{recognition problem}: to what extent does the commuting graph determine the structure it comes from? For finite groups, this question has been answered in several strong forms. Abdollahi et al.~\cite{AAM2006} conjectured that every finite non-abelian simple group is determined by its non-commuting graph (the complement of the CG), and this was confirmed by Solomon and Woldar~\cite{Solomon2013}. In a different direction,  Giudici and Parker~\cite{GiudiciParker2013} showed that the diameters of connected commuting graphs of finite groups have no universal upper bound.

For Lie algebras, the theory is still in its early stages. In 2017, Wang and Xia~\cite{Wang2017} initiated a systematic study by introducing CGs of finite-dimensional Lie algebras. They also determined the diameters, or provided upper and lower bounds for the diameters, of these graphs for simple Lie algebras over an algebraically closed field of characteristic zero. The authors of~\cite{HHVT2024} then studied CGs of a broader class of finite-dimensional Lie algebras, obtaining connectivity criteria and diameter estimates. The closed neighborhood of a vertex $x$ in the CG of a Lie algebra $\G$ is the centralizer of $x$ with the center $\Zz(\G)$ removed, and cliques correspond to families of pairwise commuting elements. Graph-theoretic statements about the CG thus translate into statements about the arrangement of abelian subalgebras of $\G$.

The recognition philosophy is native to Lie theory itself. By the Cartan--Killing classification, a semisimple Lie algebra over an algebraically closed field of characteristic zero is determined up to isomorphism by its root system (equivalently, its Cartan matrix or Dynkin diagram) and, by Serre's theorem, can be reconstructed from these data~\cite{Humphreys1972}; generalized Cartan matrices lead to Kac--Moody algebras~\cite{Kac1990}. For solvable Lie algebras, however, no combinatorial data of comparable strength are available, whereas commuting graphs are defined for every Lie algebra. This motivates the following version of the recognition problem.

\begin{quote}
    \noindent {\bf Main question.} \emph{Which classes of finite-dimensional Lie algebras, and which of their structural invariants (e.g., the dimension of the center or the existence of abelian ideals of small codimension), can be recognized from their commuting graphs?}
\end{quote}

Unfortunately, for CGs of Lie algebras we encounter an obstacle that does not arise for finite groups. Over a field $\F$ of characteristic zero, the CG $\Gamma(\G)$ of a non-abelian Lie algebra has infinitely many vertices and, as an abstract graph, forgets the linear structure of $\G$. The loss is twofold: the center fills $\Gamma(\G)$ with twin vertices (two vertices differing by a central element have the same closed neighborhood), and abstract graph isomorphisms need not respect scalar multiples, sums or dimensions of subspaces.

To remove this obstacle, we retain the linear structure inside the graph itself. We introduce \emph{linear structure graphs} (LSGs): graphs whose vertex set is a finite-dimensional $\F$-vector space $V$ without the zero vector and whose edge set is stable under the linear operations (Definition~\ref{Def2.1New}), together with \emph{LSG-monomorphisms} and \emph{LSG-isomorphisms}. 

Based on the framework of LSGs, this paper investigates two specific instances of this class of graphs associated with finite-dimensional Lie algebras over an arbitrary field $\F$ of characteristic zero (with the additional assumption of algebraic closure where necessary). These are referred to, respectively, as Lie algebra graphs (LAGs) and reduced commuting graphs (RCGs). From there, this paper will address, step by step, the Main question concerning specific classes of Lie algebras, namely, the class of commutative transitive Lie algebras (CTLAs) and the class of almost abelian Lie algebras (AALAs).

The rest of the paper is organized as follows. 
\begin{itemize}
    \item {\bf Section~\ref{Sec1New}} collects the notation used throughout. 
    
    \item {\bf Section~\ref{Sec2New}} introduces LSGs, the LSG-neighborhood of a subspace or a vertex, LSG-neighborhood sequences (Definitions~\ref{Def2.9New} and~\ref{Def2.12New}) and LSG-height (Definition~\ref{Def2.14New}). Unlike the neighborhood sequence in a CG \cite[Definition~4.1 and Theorem~4.1]{HHVT2024}, the LSG-neighborhood sequence of each vertex of an $n$-LSG $\Gamma$ terminates after at most $n-1$ steps. The first main result, Theorem~\ref{thm1}, describes the influence of the LSG-height $h_{LSG}(\Gamma)$ of $\Gamma$ on the structure and connectivity of the LSG $\Gamma$. Moreover, $h_{LSG}(\Gamma) \leq \min\{\dim(\Gamma)-1, \, \diam(\Gamma)\}$.
        
    \item {\bf Section~\ref{Sec3New}} introduces the \emph{Lie algebra graph} (LAG), a special LSG attached to a Lie algebra with trivial center. The second main result, Theorem~\ref{thm2}, characterizes LAGs of LSG-height at most $1$ as those of CTLAs (see Definition~\ref{Def3.4New}).
    
    \item {\bf Section~\ref{Sec4New}} introduces the \emph{reduced commuting graph} (RCG) $\Gamma_R(\G)$: the subgraph of the CG $\Gamma(\G)$ induced by $(\G_R)_0\co\G_R\setminus\{\boldsymbol{0}\}$, where $\G_R$ is a complementary subspace of $\Zz(\G)$ in $\G$. It removes exactly the redundancy created by the center while keeping the linear structure. The third main result, Theorem~\ref{thm3}, demonstrates that $\Gamma_R(\G)$ is well defined up to LSG-isomorphism. The fourth main result, Theorem~\ref{thm4}, shows that $\Gamma_R(\G)$ is a faithful replacement for $\Gamma(\G)$ (same connectivity, diameter and number of components). The fifth, Theorem~\ref{thm5}, shows that $\Gamma_R(\G)$ is compatible with direct-sum decompositions, so it suffices to study RCGs of indecomposable Lie algebras. Following Section~\ref{Sec3New}, Subsection~\ref{Subsec4.3New} discusses the class of AC Lie algebras (ACLAs) in Definition~\ref{Def4.6New}, which was introduced in reference~\cite{Sriw2025}. Specifically, Proposition~\ref{Prop4.8New} characterizes this class by the condition that the height of the RCG is at most 1, while also describing the RCG of ACLAs.     
    In Subsection~\ref{Subsec4.4New}, we provide a complete characterization of the RCGs of the class of AALAs (see Definition~\ref{Def4.7New}) in Theorem~\ref{thm6}, the sixth main result of the paper. 
    This result shows that the property of being almost abelian is determined by the LSG-height and LSG-dimension of the RCG.      
    
    \item {\bf Section~\ref{Sec5}}, titled ``Conclusion, open problems and declaration of AI use'', summarizes the main results of the paper; introduces Open problems~\ref{OP1New},~\ref{OP2New} and~\ref{OP3New}, and provide a declaration of AI use in Subsection~\ref{Subsec5.3New}. 
\end{itemize}

\section{Preliminaries}\label{Sec1New}
 We use the following notation throughout.
\begin{itemize}
	\item $\F$ is a field of characteristic zero and $\F_0 \co\F \setminus \{0\}$. For $n \in \N_0 \co \N\setminus\{0\}$, $\F^n$ is the $n$-dimensional $\F$-vector space and $(e_1, \ldots, e_n)$ is the standard basis of $V = \F^n$, $\F^n_0 \co \F^n \setminus \{\boldsymbol{0}\}$ ($\boldsymbol{0}$ is the zero vector).
    
    \item For any $\F$-vector space $V$, the subscript $0$ always means that the zero vector is \emph{removed}: $V_0 \co V \setminus \{\boldsymbol{0}\}$. We write $W \leq V$ if $W$ is a subspace of $V$, and $W < V$ if $W$ is a \emph{proper} subspace of $V$, i.e.\ $W \leq V$ and $W \neq V$.
    \item Suppose that $\G, \G', \mathcal{H}, \cdots$ are Lie algebras over $\F$. We write $\mathcal{H} \leq_{\text{Lie}} \G$ when $\mathcal{H}$ is a subalgebra of $\G$, and $\mathcal{H} <_{\text{Lie}} \G$ if $\mathcal{H}$ is a \emph{proper} subalgebra of $\G$, i.e.\ $\mathcal{H} \leq_{\text{Lie}} \G$ and $\mathcal{H} \neq \G$. The symbol $\mathcal{I} \triangleleft \G$ means that $\mathcal{I}$ is an ideal of $\G$. We write $f: \G \hookrightarrow_{\text{Lie}}\G'$ (resp., $f: \G \twoheadrightarrow_{\text{Lie}}\G'$ or $f: \G \overset{\cong}{\to}_{\text{Lie}}\G'$) to mean that $f$ is a Lie monomorphism (resp., Lie epimorphism or Lie isomorphism).   
	\item The symbol $\s\{x, y, \ldots\}$ is always the span of an arbitrary subset $\{x, y, \ldots\}$ in $V$. For $v \in V_0$ we have $\s\{v\} \setminus \{\boldsymbol{0}\} = \F_0 v$.
    \item Unless otherwise stated, the term ``vector space'' means a non-trivial finite-dimensional $\F$-vector space, ``Lie algebra'' means a finite-dimensional non-abelian Lie algebra over $\F$. For a given Lie algebra, when listing its Lie brackets, it is always understood that any omitted Lie brackets are trivial. $\G \oplus_{\text{Lie}} \G'$ is the direct sum of Lie algebras, and $\G \co \F^n \oplus _D \F e_0$ is the semidirect sum of $\F^n$ and $\F e_0$ induced by a non-trivial $D \in \text{Der}(\F^n) \equiv \mathfrak{gl}_n(\F)$, i.e. $e_0 \notin \F^n$, $[e_0, v] \co D(v) \in \F^n$ for all $v \in \F^n$.
\end{itemize}

\subsection{Basic notions on graphs}\label{Subsec1.1New}
Unless otherwise specified, ``graph'' means a simple, undirected graph; it may be infinite (see~\cite{BM1976,Diestel2017,West2001} for background).

\begin{definition}\label{Def1.1New}
Let $\Gamma =(V_{\Gamma}, E_{\Gamma})$ be a graph with the vertex set $V_\Gamma \neq \emptyset$ and the edge set $E_{\Gamma}$ ($\subseteq \bigl{\{}\{x, y\} \mid x, y \in V_{\Gamma}, x \neq y \bigr{\}}$). Assume that $\Sc$ is a non-empty subset of $V_{\Gamma}$.
\begin{enumerate}[(1)]
    \item\label{Def1.1-1} Two vertices $x, y$ are \emph{adjacent} if $x = y$ or $\{x, y\} \in E_{\Gamma}$. The \emph{neighborhood} of $x$ is $N(x) \co \{y \in V_{\Gamma} \mid x, y \text{ adjacent}\}$. So, $s \in N(s)$ for every vertex $s \in V_{\Gamma}$.
    
    \item\label{Def1.1-2} $\Gamma$ is called {\it complete} if every pair of vertices is adjacent.
    
    \item\label{Def1.1-3} The \emph{subgraph of $\Gamma$ induced by} $\Sc$, denoted by $\Gamma_{\Sc} = \bigl(V_{\Sc}, E_{\Sc}\bigr)$, is the graph defined as follows: $V_{\Sc} \co \Sc, \, E_{\Sc} \co \bigl{\{}\{x, y\}\in E_{\Gamma}\,|\, x, y \in \Sc\bigr{\}}$.
    
    \item\label{Def1.1-4} $\Sc$ is called a {\it clique} in $\Gamma$ if the subgraph $\Gamma_{\Sc}$ induced by $\Sc$ is complete.
    
    \item\label{Def1.1-5} A \emph{path of length $d$} ($d\in\N_0$) joining $u, v$ is a sequence $u = v_0, v_1, \cdots, v_d = v$ of vertices with $\{v_i, v_{i+1}\} \in E_{\Gamma}$, $i= 0,\dots,d-1$. Two vertices are \emph{connected} if they coincide or are joined by a path.
    
    \item\label{Def1.1-6} $\Sc$ is called a {\it connected} subset if any two of its vertices are connected. The graph $\Gamma$ is \emph{connected} if its vertex set is connected.
    
    \item\label{Def1.1-7} The maximal connected subset (in the sense of inclusion) that contains the vertex $s$ is called the {\it connected component} of $s$, denoted by $CC(s)$.
    
    \item\label{Def1.1-8} The distance function $\dist_{\Gamma} \colon V_{\Gamma} \times V_{\Gamma} \to \N \cup \{+\infty\}$ is defined as	
		\[
			\dist_{\Gamma}(x, y) \co 
			\begin{cases}
				0, & \text{if } x = y, \\
				d \geq 1, & \text{if } x \neq y \text{ and $d$ is the length of}\\
				 & \text{ a shortest path joining them}, \\
				 +\infty, & \text{if } x \neq y \text{ and they are not connected}.
			\end{cases}
		\]
        
    \item\label{Def1.1-9} The diameter $\diam(\Gamma)$ of $\Gamma$ is defined as
		\[
			\diam(\Gamma) \co 
			\begin{cases}
				\sup \{\dist_{\Gamma}(u,v) \,|\, u, v \in V_{\Gamma}\} \in \N \cup \{+\infty\}, & \text{if $\Gamma$ is connected},\\
				+\infty, & \text{otherwise}. 
			\end{cases}
		\]    
\end{enumerate}
\end{definition}

\begin{definition}\label{Def1.2New}
Let $\Gamma_i = (V_i, E_i)$ ($i=1,2$) be graphs. A map $f \colon V_1 \to V_2$ is a \emph{full homomorphism} if for all distinct $u, v \in V_1$:
\[
\{u, v\} \in E_1 \iff \bigl[f(u) = f(v)\bigr] \vee \bigl[\{f(u), f(v)\} \in E_2\bigr].
\]
An injective (resp. surjective, bijective) full homomorphism is a \emph{monomorphism} 
(resp. \emph{epimorphism}, \emph{isomorphism}). 
Any monomorphism is an isomorphism onto its image.
\end{definition}

\subsection{The commuting graph of a Lie algebra}\label{Subsec1.2New}

\begin{definition}[{\cite{HHVT2024}}]\label{Def1.3New}
Let $\G$ be a non-abelian Lie algebra with center $\Zz(\G)$. The \emph{commuting graph} (CG) $\Gamma(\G) =\big(V(\G), E(\G)\big)$ of $\G$ has vertex set $V(\G) \co \G \setminus \Zz(\G)$ and edge set $E(\G) \co \{ \{x, y\} \mid x, y \in V(\G), x \neq y, [x,y]=0 \}$.
\end{definition}

\begin{remark}\label{Re1.4New}
By~\cite[Remark 3]{HHVT2024}, if $\Gamma(\G)$ is connected then $\diam(\Gamma(\G)) \geq 2$. 
\end{remark}

\begin{proposition}[{\cite[Proposition~6]{HHVT2024}}]\label{Prop1.5New}
Let $\G$ be a decomposable Lie algebra. 
    \begin{enumerate}[(i)]
    \item\label{Pro1.5-i} If $\G = \G_1 \oplus_{\rm{Lie}} \F^k$ with $\G_1$ non-abelian, then $\Gamma(\G)$ is connected if and only if $\Gamma(\G_1)$ is, and then $\diam(\Gamma(\G)) = \diam(\Gamma(\G_1))$.
        
    \item\label{Pro1.5-ii} If $\G = \G_1 \oplus_{\rm{Lie}} \G_2$ with non-abelian Lie algebras $\G_1,\G_2$, then $\Gamma(\G)$ is connected and $2 \leq \diam(\Gamma(\G)) \leq 3$. \hfill $\square$
    \end{enumerate} 
\end{proposition}

\section{Linear structure graphs and their characteristics}\label{Sec2New}

\subsection{Linear structure graphs}\label{Subsec2.1New}

\begin{definition}[{\bf Linear structure graph}]\label{Def2.1New}
Let $\Gamma = (V_{\Gamma}, E_{\Gamma})$ be a graph for which there is a vector space $V$ with $V_{\Gamma} = V_0$. Assume that $E_{\Gamma}$ is \emph{stable under linear operations}: for all $u_1, u_2, u, v, w \in V_{\Gamma}$ and $\alpha, \beta \in \F_0$,
    \begin{enumerate}[($L_1$)]    
    \item\label{Def2.1-$L_1$} if $\alpha\neq\beta$, then $\{\alpha u,\beta u\}\in E_\Gamma$;
    
    \item\label{Def2.1-$L_2$} if $\{v,w\}\in E_\Gamma$ and $\alpha v\neq\beta w$, then $\{\alpha v,\beta w\}\in E_\Gamma$;
    
    \item\label{Def2.1-$L_3$} if $\{u_1,v\},\{u_2,v\}\in E_\Gamma$ and $\boldsymbol{0}\neq u_1+u_2\neq v$, then $\{u_1+u_2,v\}\in E_\Gamma$.
    \end{enumerate}     
    Then $\Gamma$ is called a \emph{linear structure graph} (LSG) \emph{constructed on} $V$. If $V = \F^n$ it is an \emph{$n$-LSG}, and $\dim (\Gamma) = n$.
\end{definition}

In contrast to the CGs of Lie algebras, for the LSGs, the key point is to consider LSG-adjacency across the entire space $V$, including the zero vector.

\begin{definition}[{\bf LSG-adjacency, LSG-neighborhood}]\label{Def2.2New}
Let $\Gamma$ be an LSG constructed on $V$. For $u,v\in V$ (possibly zero vectors), $u,v$ are \emph{LSG-adjacent}, written $u\sim _{LSG}v$, if $\boldsymbol{0} \in \{u, v\}$ or $u = v$, or $\{u,v\}\in E_{\Gamma}$. For $v\in V_0$, the \emph{LSG-neighborhood} is $N_{LSG}(v) \co \{u\in V\mid u\sim _{LSG}v\}$.
\end{definition}

The LSG-adjacency is reflexive and symmetric but not transitive; most importantly, however, it possesses linear stability.

\begin{proposition}[{\bf Linear stability of the LSG-adjacency}]\label{Prop2.3New} 
Let $\Gamma$ be an LSG constructed on $V$. Then the following statements hold.
\begin{enumerate}[(i)]
    \item\label{Pro2.3-i} For every $u, u_1, u_2, v \in V$ (possibly zero vector) and any $\alpha, \beta \in \F$, we have
    \[(u \sim_{LSG} v) \Rightarrow (\alpha u \sim_{LSG} \beta v); \, \, \bigl[(u_1\sim_{LSG}v) \wedge (u_2\sim_{LSG}v)\bigr] \Rightarrow (u_1+u_2\sim_{LSG}v).\]
    
    \item\label{Pro2.3-ii} For every $u, v \in V_0$ and any $\alpha \in \F_0$, we have
    \[N_{LSG}(\alpha v) = N_{LSG}(v); \, \, \bigl(u \in N_{LSG}(v)\bigr) \iff \bigl(v \in N_{LSG}(u)\bigr).\]
    
    \item\label{Pro2.3-iii} For every $v \in V_0$, $N_{LSG}(v) = \{\boldsymbol{0}\}\sqcup N(v)$ is a subspace of $V$ containing $\F v$
    (here $N(v)$ is the ordinary neighborhood of $v$ in $\Gamma$).
\end{enumerate} 
\end{proposition}

\begin{proof}  
Assertion (i) follows directly from Definitions~\ref{Def2.1New} and~\ref{Def2.2New} after carefully considering all possibilities for the vectors $u, u_1, u_2, v \in V$ and scalars $\alpha, \beta \in \F$. We only need to prove the remaining two assertions.
\begin{itemize}
    \item {\bf Proof of~\eqref{Pro2.3-ii}.} By Assertion (i), for any $u, v \in V_0$ and any $\alpha \in \F_0$, we have
    \[(u\sim_{LSG} v) \Rightarrow (u\sim_{LSG}\alpha v) \Rightarrow \bigl(u\sim_{LSG}\alpha^{-1}(\alpha v)) \Rightarrow (u \sim_{LSG} v).\]
    Hence, $(u\sim_{LSG} v) \Leftrightarrow (u \sim_{LSG} \alpha v)$, i.e. $N_{LSG}(\alpha v)=N_{LSG}(v)$. Because of the symmetry of the LSG-adjacency in Definition~\ref{Def2.2New}, we also have 
    \[\bigl(u \in N_{LSG}(v)\bigr) \iff \bigl(v \in N_{LSG}(u)\bigr).\]
        
    \item {\bf Proof of~\eqref{Pro2.3-iii}.} By Assertion~\eqref{Pro2.3-i}, for every $v \in V_0$, $N_{LSG}(v)$ is closed under addition and scalar multiplication and it contains $\boldsymbol{0}$ and $v$. Hence, it is a subspace of $V$ containing $\F v$. Finally, $N_{LSG}(v)\setminus\{\boldsymbol{0}\}=N(v)$ by Definitions~\ref{Def1.1New}\eqref{Def1.1-1} and~\ref{Def2.2New}. 
\end{itemize} The proposition is proven. 
\end{proof}

\begin{definition}[{\bf Linear structure subgraph}]\label{Def2.4New}
Let $\Gamma$ be an LSG constructed on $V$ and $0 \ne W \leq V$. Then the subgraph $\Gamma_{W}$ induced by $W_0$ is an LSG constructed on $W$ and it is called a \emph{linear structure subgraph} (LSsubG) of $\Gamma$. If $\Gamma_W$ is connected, we call it \emph{LSsubG-connected} and $W$ is called an {\it LSG-connected subspace} of $V$. If $W_0$ is also a connected component of $V_0$, then $W$ is called an \emph{LSG-connected component} in $V$ and $\Gamma_W$ is called an {\it LSsubG-connected component} of the LSG $\Gamma$.
\end{definition}

\begin{definition}[{\bf Direct sum of LSGs}]\label{Def2.5New}
Let $\Gamma_i$ be LSGs constructed on $V_i$ ($i=1,2$) and $V := V_1 \oplus V_2$. The graph $\Gamma = (V_{\Gamma}, E_{\Gamma})$ with
\[V_{\Gamma} = V_0, \, E_{\Gamma} = \Bigl{\{} \{u, v\} \subset V_0\,|\, u=(u_1,u_2) \ne v=(v_1,v_2), u_1 \sim_{LSG} v_1, u_2 \sim_{LSG} v_2 \Bigr{\}}\]
is called the {\it LSG-direct sum} of LSGs $\Gamma_1$ and $\Gamma_2$, denoted by $\Gamma_1 {\oplus}_{LSG} \Gamma_2$.    
\end{definition}

\begin{proposition}\label{Prop2.6New}
$\Gamma_1 \oplus_{LSG} \Gamma_2$ is an LSG constructed on $V_1 \oplus V_2$; in particular $\dim(\Gamma_1 \oplus_{LSG} \Gamma_2) = \dim(\Gamma_1) + \dim(\Gamma_2)$.
\end{proposition}

\begin{proof}
Clearly, for $(u_1, u_2) = u \neq v = (v_1, v_2)$, 
\[\{u, v\} \in E_{\Gamma} \iff \bigl(u_1 \sim_{LSG} v_1\bigr) \wedge \bigl(u_2 \sim_{LSG} v_2\bigr).\]   
By Proposition~\ref{Prop2.3New}(i), both of the LSG-adjacencies in $\Gamma_1, \, \Gamma_2$ are linearly stable. Therefore, $(L_1)$--$(L_3)$ are satisfied for the adjacency of $\Gamma = \Gamma_1 \oplus_{LSG} \Gamma_2$. In other words, $\Gamma$ is an LSG constructed on $V = V_1 \oplus V_2$, and 
$\dim(\Gamma_1 \oplus_{LSG} \Gamma_2) = \dim(\Gamma_1) + \dim(\Gamma_2)$.
\end{proof}

\begin{definition}[{\bf LSG-monomorphisms, LSG-isomorphisms}]\label{Def2.7New}
Let $\Gamma, \Gamma'$ be LSGs constructed on $V, V'$, respectively, and $f\colon V \hookrightarrow V'$ be a linear monomorphism. 
\begin{enumerate}[(1)]
    \item\label{Def2.7-1} {\bf LSG-monomorphisms}: If $\bar{f}\co f|_{V_0}$ is also a graph monomorphism $\Gamma \to \Gamma'$ (Definition~\ref{Def1.2New}), then it is called the \emph{LSG-monomorphism induced by $f$}, written $\bar{f} = LSG(f)\colon \Gamma \hookrightarrow _{LSG} \Gamma'$. In particular, the natural embeddings $j_i \colon V_i \hookrightarrow V_1\oplus V_2$ induce \emph{canonical LSG-embeddings} $LSG(j_i) \colon\Gamma_i \hookrightarrow_{LSG} \Gamma_1\oplus_{LSG}\Gamma_2$, ($i = 1, 2$).
    
    \item\label{Def2.7-2} {\bf LSG-isomorphisms}: If $f\colon V \overset{\cong}{\to} V'$ is also a linear isomorphism, then $LSG(f)$ is an \emph{LSG-isomorphism}, written $LSG(f) \colon \Gamma \overset{\cong}{\to}_{LSG} \Gamma'$.

    \item\label{Def2.7-3} {\bf LSG-inavariant}: A property (P) of $\Gamma$ is called an {\it LSG-invariant} of $\Gamma$ if every $\Gamma'$ that is LSG-isomorphic to $\Gamma$ also possesses that property.
\end{enumerate}  
\end{definition}

\begin{example}[{\bf LSG-complete}]\label{Ex2.8New}
Clearly, for any vector space $V$, the complete graph $\boldsymbol{K}_{V_0}$ is an LSG constructed on $V$ and it is called the {\it LSG-complete on $V$}. When $V=\F^n$, we rewrite $\boldsymbol{K}_{V_0}$ as $\boldsymbol{K}_{\F^n_0}$. 
\end{example}

\subsection{LSG-closed subspaces and the LSG-neighborhood sequences}\label{Subsec2.2New}

\begin{definition}[{\bf LSG-neighborhood of a subspace, LSG-closed subspace}]\label{Def2.9New} 
Let $\Gamma = (V_{\Gamma}, E_{\Gamma})$ be an LSG constructed on $V$ and $S_V \co \{W\,|\, W \leq V\}$. 
    \begin{enumerate}[(1)]
    \item\label{Def2.9-1} {\bf The LSG-neighborhood map and LSG-neighborhood of a subspace}. Consider the map $\Nn_{LSG} \colon S_V \to S_V$ defined by
    \begin{equation*}
    \Nn_{LSG} (W) \coloneqq \begin{cases}
    \{\boldsymbol{0}\} & \text{if} \, \, W = \{\boldsymbol{0}\};\\
    \sum_{w\in W_0} N_{LSG}(w) & \text{if} \, \, W \ne \{\boldsymbol{0}\}.
    \end{cases}            
    \end{equation*} This map is called the {\it LSG-neighborhood map} of $\Gamma$. The image $\Nn_{LSG} (W)$ is called the {\it LSG-neighborhood} of the subspace $W$ in $V$.
        
    \item\label{Def2.9-2} {\bf LSG-neighborhood of a vertex}. For any vertex $s \in V_{\Gamma}$, the LSG-neighborhood $\Nn_{LSG}(\F s)$ of the subspace $\F s$ is called the {\it LSG-neighborhood of $s$}, denoted by $\Nn_{LSG}(s)$.
        
    \item\label{Def2.9-3} {\bf LSG-closed subspace, LSG-closed vertex and the family of all LSG-closed subspaces} 
        \begin{itemize}
            \item {\bf an LSG-closed subspace}: If $W = \Nn_{LSG}(W)$ then it is called an {\it LSG-closed subspace} of $V$.
            
            \item {\bf an LSG-closed vertex}: For any vertex $s \in V_{\Gamma}$, if $\F s$ is an LSG-closed subspace of $V$, i.e. $\Nn_{LSG}(s) = \F s$, then $s$ is called an {\it LSG-closed vertex} in $V_{\Gamma}$.
            
            \item {\bf Family of all LSG-closed subspaces of $V$}: 
            \[CS_{LSG}(V) \co \{W \leq V \mid W \, \, \text{is an LSG-closed subspace of}\, \, V\}.
            \]            
        \end{itemize}               
    \end{enumerate}
\end{definition}

\begin{remark}\label{Re2.10New} For any LSG $\Gamma$, by Proposition~\ref{Prop2.3New} and Definition~\ref{Def2.9New}, we have
\begin{enumerate}[(1)]
    \item\label{Re2.10-1} $\Nn_{LSG}(s) = N_{LSG}(v) = \{\boldsymbol{0}\}\sqcup N(v); \, \forall v \in \F_0 s$ for any vertex $s$ of $\Gamma$. 
    
    \item\label{Re2.10-2} A subspace of an LSG-closed subspace need not be LSG-closed.
\end{enumerate}
\end{remark}

\begin{proposition}\label{Prop2.11New}
Let $\Gamma = (V_{\Gamma}, E_{\Gamma})$ be an LSG constructed on a vector space $V$, $s \in V_{\Gamma}$ and $W \leq V$. Then the following statements hold.
    \begin{enumerate}[(i)]
    \item\label{Pro2.11-i} $W \leq \Nn_{LSG}(W)$ and the map $\Nn_{LSG}$ is monotonically increasing: 
    \[(W \leq W') \Longrightarrow \bigl(\Nn_{LSG}(W) \leq \Nn_{LSG}(W')\bigr).\]
        
    \item\label{Pro2.11-ii} Regarding LSG-closed subspaces of $\Gamma$ and the family $CS_{LSG}(V)$, we have
    \begin{enumerate}[(a)]
        \item $\{0\}, V \in CS_{LSG}(V)$.
        \item $(s \, \, \text{is LSG-closed}) \, \, \Longleftrightarrow \bigl(N(v) = \F_0 s; \, \forall v\in\F_0 s\bigr)$.
        \item Any intersection of LSG-closed subspaces is again LSG-closed.
    \end{enumerate}               
    \end{enumerate}    
\end{proposition}

\begin{proof} We will show that each assertion of the proposition is true.
\begin{enumerate}[(i)]
    \item {\bf Proof of~\eqref{Pro2.11-i}.} First of all, Assertion (i) is obvious by Definition~\ref{Def2.9New}\eqref{Def2.9-1}.
    
    \item {\bf Proof of~\eqref{Pro2.11-ii}.} Next, for Assertion (ii), it is clear that
    \begin{enumerate}[(a)]
        \item $\Nn_{LSG}(\{\boldsymbol{0}\}) = \{\boldsymbol{0}\}$ and $\Nn_{LSG}(V) = V$, so $\{\boldsymbol{0}\}, V$ are themselves LSG-closed subspaces. This means (a) is correct.
        \item By Remark~\ref{Re2.10New}, $\Nn_{LSG}(s) = N_{LSG}(v) = \{\boldsymbol{0}\}\sqcup N(v); \, \forall v \in \F_0 s, \, \forall s \in V_{\Gamma}$. Therefore, $(s \, \, \text{is LSG-closed}) \, \, \Longleftrightarrow \bigl(N(v) = \F_0 s; \, \forall v\in\F_0 s\bigr)$; i.e. (b) is correct. 
        \item Let $\{W_i\}_{i\in I} \subseteq CS_{LSG}(V)$ and $W=\bigcap_{i\in I} W_i$, $I$ is an arbitrary index set. Since $W \leq \Nn_{LSG}(W) \co \sum_{w \in W_0} N_{LSG}(w) \leq \Nn_{LSG}(W_i)= W_i$ for all $i$, so $W \le \Nn_{LSG}(W) \leq \bigcap_{i\in I} W_i = W$. That means (c) holds.
    \end{enumerate}
\end{enumerate} Thus, Proposition~\ref{Prop2.11New} is proven. 
\end{proof}

\begin{definition}[{\bf LSG-neighborhood sequence}]\label{Def2.12New}
Let $\Gamma$ be an LSG constructed on $V$ and $0\neq W\leq V$. Put $\Nn^0_{LSG}(W) \coloneqq W$ and $\Nn^{k+1} _{LSG}(W) \coloneqq \Nn_{LSG}\bigl(\Nn^{k}_{LSG}(W)\bigr)$. The sequence $(\Nn^k_{LSG}(W))_{k\in\N}$ is called the \emph{LSG-neighborhood sequence of $W$}. In particular, for $W=\F s$ we rewrite $(\Nn^k_{LSG}(\F s))_{k\in\N}$ as $(\Nn^k_{LSG}(s))_{k\in\N}$ and call it the {\it LSG-neighborhood sequence of $s \in V_{\Gamma}$}.
\end{definition}

\begin{proposition}[{\bf Strict increase and termination of the LSG-neighborhood sequence}]\label{Prop2.13New}
Let $\Gamma$ be an $n$-LSG constructed on $\F^n$ ($n \geq 1$).
\begin{enumerate}[(i)]
    \item\label{Pro2.13-i} For any non-trivial $W \leq \F^n$, 
    there is a natural number $h = h(W) < n$  such that the LSG-neighborhood sequence of $W$ is strictly increasing and terminates after $h(W)$ steps: 
        \begin{equation}\label{Ser2.1New}
        \Nn^0_{LSG}(W)<\cdots< \Nn^{h(W)}_{LSG}(W) = \Nn^{h(W)+1}_{LSG}(W)
        \end{equation}
    \item[] Furthermore, $\Nn^{h(W)}_{LSG}(W)$ is the smallest LSG-closed subspace (under the inclusion relation) that contains $W$.
    \item\label{Pro2.13-ii} When $W = \F s$ for any vertex $s$ in $\Gamma$, there is a natural number $h = h(s) < n$ such that the LSG-neighborhood sequence of $s$ is strictly increasing and terminates after $h(s)$ steps:
        \begin{equation}\label{Ser2.2New}
        \Nn^0_{LSG}(s)< \cdots < \Nn^{h(s)}_{LSG}(s)= \Nn^{h(s)+1}_{LSG}(s)
        \end{equation}
    \item[] Furthermore, $\Nn^{h(s)}_{LSG}(s)$ is the smallest LSG-closed subspace (under the inclusion relation) that contains $s$.
\end{enumerate}
\end{proposition}

\begin{proof} It is clear that Statement~\eqref{Pro2.13-ii} follows from Statement~\eqref{Pro2.13-i} by applying it with $W = \F s$. Therefore, we only need to prove Statement~\eqref{Pro2.13-i}.
\begin{itemize}
    \item {\bf Proof of Relation~\eqref{Ser2.1New}.} By Proposition~\ref{Prop2.11New}, $\Nn^k_{LSG}(W) \leq \Nn^{k+1}_{LSG}(W)$, and when two consecutive terms are equal, the sequence terminates from that point on. Let $h = h(W)$ be the least index with $\Nn_{LSG}^h (W) =\Nn_{LSG}^{h+1}(W)$. The terms $\Nn_{LSG}^0(W),\dots,\Nn_{LSG}^h(W)$ are strictly increasing non-zero subspaces of $\F^n$, so their dimensions are strictly increasing in $[1,n]$, whence $h+1\leq n$, i.e. $h = h(W) < n = \dim(\Gamma)$ for every non-trivial subspace $W$ of $\F^n$. 
    \item {\bf Prove $\Nn^{h(W)}_{LSG}(W)$ is the smallest LSG-closed subspace containing $W$}\\ 
    By Definition~\ref{Def2.9New}\eqref{Def2.9-3},  $\Nn^{h(W)}_{LSG}(W)$ is an LSG-closed subspace of $\F^n$ containing $W$. For any $W' \in CS_{LSG}(\F^n)$, if $W \leq W'$ then Proposition~\ref{Prop2.11New}\eqref{Pro2.11-i} gives $\Nn^k_{LSG}(W) \leq \Nn^k_{LSG}(W') = W'$ for all $k \geq 1$. In particular, $\Nn^{h(W)}_{LSG}(W) \leq W'$. In other words, $\Nn^{h(W)}_{LSG}(W)$ is the smallest LSG-closed subspace containing $W$.
\end{itemize} The proposition is proven.
\end{proof}

\begin{definition}[{\bf LSG-height, smallest LSG-closed subspace}]\label{Def2.14New}
Let $\Gamma$ be an $n$-LSG on $\F^n$ ($n \geq 1$).
    \begin{enumerate}[(1)]
        \item\label{Def2.14-1} {\bf The LSG-height of $W$, the smallest LSG-closed subspace containing $W$}. For any $W \leq \F^n$ ($W \ne \{\boldsymbol{0}\}$), we define
        \begin{itemize}
            \item {\bf The LSG-height of $W$}: $h = h(W) (\in \N)$ in sequence~\eqref{Ser2.1New} is called the \emph{height} of $W$, denoted by $h_{LSG}(W)$. Of course, $0 \leq h_{LSG}(W) < n = \dim(\Gamma)$.
            \item {\bf The smallest LSG-closed subspace containing} $\boldsymbol{W}$. The subspace 
            $\Nn^{h_{LSG}(W)}_{LSG}(W)$ is called 
            the {\it smallest LSG-closed subspace containing} $W$, denoted by $\s_{LSG}(W)$. 
        \end{itemize}
        
        \item\label{Def2.14-2} {\bf The LSG-height of a vertex, the smallest LSG-closed subspace containing a vertex, the LSG-closed dimension of a vertex}. When $W=\F s$ for a vertex $s\in V_{\Gamma}$, we consider the sequence~\eqref{Ser2.2New}. 
        \begin{itemize}
            \item {\bf The LSG-height of $s$}: The LSG-height $h_{LSG}(\F s)$ of $W = \F s$ is called the {\it LSG-height} of $s$, denoted by $h_{LSG}(s)$. Of course, $0 \leq h_{LSG}(s) < n = \dim(\Gamma)$.
            \item {\bf The smallest LSG-closed subspace containing $s$}: $\s_{LSG}(\F s)$ is also called the {\it smallest LSG-closed subspace containing} $s$, denoted by $\s_{LSG}(s)$. 
            \item {\bf The LSG-closed dimension of $s$}: The dimension $\dim \s_{LSG}(s)$ is also called the {\it LSG-closed dimension} of $s$, denoted by $\dim_{LSG}(s)$.            
        \end{itemize}
        
        \item\label{Def2.14-3} {\bf The LSG-height of $\Gamma$}. The {\it LSG-height of $\Gamma$}, denoted by $h_{LSG}(\Gamma)$, is defined as follows: $h_{LSG}(\Gamma) \co \max\{h_{LSG}(s) \, |\, s \in \F^n_0\}$. Of course, $0 \leq h_{LSG}(\Gamma) < n$. 
         
        \item\label{Def2.14-4} {\bf The LSG-closed dimension of $\Gamma$}. The {\it LSG-closed dimension of $\Gamma$}, denoted by $\dim_{LSG}(\Gamma)$, is defined as follows: $\dim_{LSG}(\Gamma) \co \max\{\dim_{LSG}(s) \mid s \in \F^n_0\}$. Of course, $1 \leq \dim_{LSG}(\Gamma) \leq n$.        
    \end{enumerate}
\end{definition}

\begin{remark}\label{Re2.15New} We have the following three observations, which are obvious but important.
\begin{enumerate}[(1)]
    \item\label{Re2.15-1} The terminology ``LSG-closed subspace'' and related terms follows the traditional naming convention for closure operators and closed elements in ordered sets: by Propositions~\ref{Prop2.11New} and~\ref{Prop2.13New}, the map $W \mapsto \s_{LSG}(W)$ is a closure operator on the lattice of subspaces of $V$, and its closed elements are exactly the LSG-closed subspaces.
    
    \item\label{Re2.15-2} By Proposition~\ref{Prop2.13New}, Definition~\ref{Def2.14New}, for any vertex $s$, we have
\[\bigl(N_{LSG}(s) = \F s\bigr) \Leftrightarrow \bigl(\dim_{LSG}(s) = 1\bigr); \, \, \bigl(\F s < N_{LSG}(s)\bigr) \Leftrightarrow \bigl(1 < \dim_{LSG}(s) \leq n\bigr).\]
    \item\label{Re2.15-3} For any LSG $\Gamma$ constructed on a vector space $V$, the LSG-height and the LSG-closed dimension of $\Gamma$ are LSG-invariants of $\Gamma$.
\end{enumerate}
\end{remark}

\begin{notation}\label{No2.16New}
Let $\widetilde{s} \co \F_0 s$ denote the equivalence class of $s\in\F_0^n$ ($1 \leq n \in \N$) under the proportionality relation $\sim_{\F}$: $\forall u, v \in \F^n_0, (u\sim_{\F}v) \Leftrightarrow \bigl(\exists \alpha \in \F_0: v = \alpha u\bigr)$. The symbol $\mathbb{P}^n_{\F} \co \F^{n+1}_0/\sim_{\F} \, = \{ \widetilde{s} \mid s \in \F_0^{n+1}\}$ is the $n$-dimensional $\F$-projective space. 
We always have $\F_0^{n} \equiv \F_0^{n}\times \{0\} \subsetneq \F_0^{n+1}$, $\mathbb{P}_{\F}^{n-1} \subsetneq \mathbb{P}_{\F}^n$ and 
$(\F^{n+1}_0\setminus\F^n_0)/\sim_{\F} \, \equiv \mathbb{P}_{\F}^n \setminus \mathbb{P}_{\F}^{n-1}$.
\end{notation}

\begin{theorem}\label{thm1} 
Let $\Gamma = \bigl(V_{\Gamma}, E_{\Gamma}\bigr)$ be an $n$-LSG constructed on $\F^n$ ($n \geq 1$).
	\begin{enumerate}[(i)]
		\item\label{thm1-i} {\bf The connected component of vertex $s$ is contained in $\boldsymbol{\s_{LSG}(s)}$:} For every $s \in \F^n_0$, the connected component of $s$ satisfies
		\begin{equation}\label{For2.3New}
			CC(s) \subseteq \s_{LSG}(s)\setminus\{\boldsymbol{0}\}.
		\end{equation}
        
		\item\label{thm1-ii} {\bf The influence of the LSG-height $\boldsymbol{h_{LSG}(\Gamma)}$ on $\boldsymbol{\Gamma}$:} The structure of $\Gamma$ is determined by its $h_{LSG}(\Gamma)$ in the following sense.
		\begin{enumerate}[(a)]
			\item Assume that $\boldsymbol{h_{LSG}(\Gamma) = 0.}$ Then $CC(s)=\F_0 s$ for every $s$. Hence
			\begin{equation}\label{For2.4New}
				\Gamma \cong_{LSG} \bigsqcup_{\widetilde{s} \in \F^n_0/\sim_{\F} \, = \,  \mathbb{P}_{\F}^{n-1}} \boldsymbol{K}_{\widetilde{s}},
			\end{equation}			
            ($\Gamma$ is connected) $\iff \, \, (n = 1) \wedge \bigl(\Gamma \cong_{LSG}  \boldsymbol{K}_{\F_0}\bigr)$.
            
			\item Assume that $\boldsymbol{h_{LSG}(\Gamma) = 1.}$ Then $n \geq 2$. Furthermore, every connected component of $\Gamma$ is an LSG-complete $\boldsymbol{K}_{W_0}$ with $W = N_{LSG}(s)$ for some $s \in V_{\Gamma}$, $2 \le \dim_{LSG}(\Gamma) \le n$ and ($\Gamma$ is connected) $\iff \, \, \dim_{LSG}(\Gamma) = n$. 
            More specifically,
			\begin{itemize}
				\item {\bf Case 1:} $\boldsymbol{\dim_{LSG}(\Gamma) = n \geq 2}$. Then $\Gamma = \boldsymbol{K}_{\F^n_0}$ and $\diam(\Gamma) = 1$. 
				\item {\bf Case 2:} $\boldsymbol{\dim_{LSG}(\Gamma) = n-1 \geq 2}$. Then $n \geq 3$ and $\Gamma$ is disconnected.    
                Furthermore, $\boldsymbol{K}_{\F^{n-1}_0}$ is the unique LSG-connected component of dimension $n-1$ in $\Gamma$, and all other components are $1$-dimensional. In that case, 
                \begin{equation}\label{For2.5New}
					\Gamma \cong_{LSG} \boldsymbol{K}_{\F_0^{n-1}} \sqcup \Bigl( \bigsqcup_{\widetilde{s} \in (\F^n_0 \setminus \F^{n-1}_0)/\sim_{\F} \, = \, \mathbb{P}^{n-1}_{\F} \setminus \mathbb{P}^{n-2}_{\F}} \boldsymbol{K}_{\widetilde{s}}\Bigr)
				\end{equation}
                \item {\bf Case 3:} $\boldsymbol{2 \leq \dim_{LSG}(\Gamma) \leq n-2}$. Then $n \ge 4$ and $\Gamma$ is disconnected. Moreover, if $2 \dim_{LSG}(\Gamma) > n$, then $\boldsymbol{K}_{\F^{\dim_{LSG}(\Gamma)}_0}$ is the unique LSG-connected component of dimension $\dim_{LSG}(\Gamma)$ in $\Gamma$. 	
				\end{itemize}
                
            \item Assume that $\boldsymbol{2 \leq h_{LSG}(\Gamma) < n.}$ Then $n \geq 3$ and $\Gamma$ is not complete. Furthermore,
				\begin{equation}\label{For2.6New}
				(\text{$\Gamma$ connected}) \Rightarrow \bigl(\s_{LSG}(s) = \F^n, \, \, \forall s \in \F^n_0\bigr) \wedge \bigl(\diam(\Gamma) \geq 2\bigr)
				\end{equation}
		\end{enumerate}
        
        \item\label{thm1-iii} {\bf LSG-height and diameter of $\boldsymbol{\Gamma}$:} $0\leq h_{LSG}(\Gamma) \leq \min\bigl{\{}n-1, \diam(\Gamma)\bigr{\}}$. Moreover,
        $\bigl(\Gamma \, \, \text{is connected, not complete} \bigr) \Longrightarrow \bigl[ \bigl(2 \leq h_{LSG}(\Gamma)\bigr) \wedge (n \geq 3)\bigr]$.        
	\end{enumerate}
\end{theorem}

\begin{proof}[{\bf Proof of Theorem~\ref{thm1}}] We will prove each assertion in turn.
\begin{itemize}
    \item {\bf Proof of~\eqref{thm1-i}.} Let $u\in CC(s)\setminus\{s\}$ and  $s=u_0,u_1,\dots,u_h=u$ be a path in $CC(s)$. Then $u_i\in N_{LSG}(u_{i-1})$, so inductively $u_i\in\Nn^i_{LSG}(s)$. By Proposition~\ref{Prop2.13New}, $\Nn^i_{LSG}(s)\leq \Nn^{h(s)}_{LSG}(s)=\s_{LSG}(s)$ for all $i$. Hence $u\in\s_{LSG}(s)$ and (i) holds.
    
    \item {\bf Proof of~\eqref{thm1-ii}}
    \begin{itemize}
        \item {\bf (ii)(a)}. Assume that $h_{LSG}(\Gamma) = 0$, i.e. $h_{LSG}(s)=0$ and $\F s = N_{LSG}(s)$ for all $s \in V_{\Gamma}$. In other words, \ $s \sim_{LSG} v$ if and only if $v \in \F_0 s$. Thus $\F_0 s$ is a clique and a connected component of $\Gamma$. Furthermore, the components are in bijection with $\mathbb{P}^{n-1}_\F$. This gives Relation~\eqref{For2.4New}. Moreover, connectedness forces a single component, i.e. $n=1$. That means (ii)(a) holds.
        \item {\bf (ii)(b)}. Assume that $h_{LSG}(\Gamma)=1$. Now $h_{LSG}(s)\le1$ for all $s$. So $\Nn^2_{LSG}(s)=\Nn^1_{LSG}(s)$, i.e.\ $N_{LSG}(w)\leq N_{LSG}(s)$ for all $w\in N_{LSG}(s)_0$. Applying the same argument to the vertex $w$ gives the reverse inclusion, so
        \begin{equation}\label{For2.7New}
        N_{LSG}(w)=N_{LSG}(s); \, \, \forall w\in N_{LSG}(s)_0 
        \end{equation}
        Hence, any two distinct $u,w\in N_{LSG}(s)_0$ are adjacent and $N_{LSG}(s)_0$ is a clique containing $s$. So $N_{LSG}(s)_0\subseteq CC(s)$. Together with Relation~\eqref{For2.3New} and $\s_{LSG}(s)=\Nn^1_{LSG}(s)=N_{LSG}(s)$, this gives $CC(s)=N_{LSG}(s)_0$, i.e. an LSsubG-complete $\boldsymbol{K}_{W_0}$ with $W=N_{LSG}(s)$. By Relation~\eqref{For2.7New}, two such subspaces coincide or meet only in $\{\boldsymbol{0}\}$, i.e. they form a vector space partition of $\F^n$. Since $h_{LSG}(\Gamma)=1$, there is $s_0$ with $h_{LSG}(s_0)=1$ and $\dim N_{LSG}(s_0)=\dim_{LSG}(\Gamma)\ge 2$ (Remark~\ref{Re2.15New}), so $n\ge2$. 
        \begin{itemize}  
            \item[$\ast$] {\bf Case 1:} If $\Gamma$ is connected then $CC(s_0)=\F_0^n$, so $\dim_{LSG}(\Gamma)=n$. Conversely $\dim_{LSG}(\Gamma)=n$ gives $N_{LSG}(s_0)=\F^n$ and $\Gamma=\boldsymbol{K}_{\F^n_0}$ is complete with diameter $1$.
            \item[$\ast$] {\bf Case 2:} Suppose $\dim_{LSG}(\Gamma)=n-1\ge2$, so $n\ge 3$. Take $s_0 \in V_{\Gamma}$ with $W=N_{LSG}(s_0)\cong \F^{n-1}$. For any $u\notin W$, $CC(u) \ne CC(s_0)$,        so $N_{LSG}(u)\cap W=\{\boldsymbol{0}\}$ and $\dim N_{LSG}(u)\le1$, i.e.\ $CC(u)=\widetilde u$. The $1$-dimensional components are indexed by the lines not contained in $W$, i.e. by $(\F^n _0 \setminus \F^{n-1}_0)/\sim \, \equiv \, \mathbb{P}^{n-1}_{\F} \setminus \mathbb{P}^{n-2}_{\F}$; this is Relation~\eqref{For2.5New}.
            \item[$\ast$] {\bf Case 3:} Next, if $2\le k\co\dim_{LSG}(\Gamma)\le n-2$, then $\Gamma$ is disconnected and every other component $W'$ meets $\F^k$ trivially, so $\dim W'\le n-k$. Therefore, if $2k>n$, then $\dim W' < k$, which gives uniqueness of the LSG-connected component of dimension $k = \dim_{LSG}(\Gamma)$ in $\Gamma$. 
        \end{itemize} That means (ii)(b) holds.      
        \item {\bf (ii)(c)}. Assume that $2 \leq h_{LSG}(\Gamma) < n$. Then, of course, $n \geq 3$. If $\Gamma$ were complete, then $N_{LSG}(s)=\F^n$ for all $s$, so $h_{LSG}(s)\le1$, contradicting $h_{LSG}(\Gamma)\ge2$. That means $\Gamma$ is not LSG-complete. By Relation~\eqref{For2.3New}, $CC(s) \subseteq \s_{LSG}(s)\setminus\{\boldsymbol{0}\}$. If $\Gamma$ is connected then $CC(s) = \F^n_0 \subseteq \s_{LSG}(s)\setminus\{\boldsymbol{0}\}$. From that, there must be $\s_{LSG}(s)=\F^n$.
        Moreover, since $\Gamma$ is connected but not complete, there are two distinct non-adjacent vertices $u, v$, so $\dist_{\Gamma}(u, v) \geq 2$ and $\diam(\Gamma) \geq 2$. This is Relation~\eqref{For2.6New}, so (ii)(c) also holds.
    \end{itemize}
    
    \item {\bf Proof of~\eqref{thm1-iii}.} By Definition~\ref{Def2.14New}\eqref{Def2.14-3}, we have $0\leq h_{LSG}(\Gamma) \leq n-1$. Furthermore, when $\diam(\Gamma) = +\infty$, there is nothing further to prove. Thus, we need only consider the case where $\Gamma$ is connected with $\diam(\Gamma) \in \N$. We will prove that $h_{LSG}(s) \leq \diam(\Gamma)$ for any vertex $s$. Fix an arbitrary vertex $s \in \F^n_0$.
    \begin{itemize}
        \item[$\ast$] If $h_{LSG}(s)=0$, there is nothing to prove.
        
        \item[$\ast$] Assume that $h_{LSG}(s)>0$. The sequence~\eqref{Ser2.2New} of LSG-neighborhoods of $s$ (see Proposition~\ref{Prop2.13New}) is strictly increasing up to step $h_{LSG}(s)$, so there is a vertex $u_0 \in \Nn^{h_{LSG}(s)}_{LSG}(s) \setminus \Nn^{h_{LSG}(s)-1}_{LSG}(s)$. Every vertex joined to $s$ by a path of length $k$ lies in $\Nn^{k}_{LSG}(s)$, as shown in the proof of Item~(i), so no path of length less than $h_{LSG}(s)$ can connect $s$ to $u_0$. Therefore, by Definition~\ref{Def1.1New}\eqref{Def1.1-8},\eqref{Def1.1-9}, $h_{LSG}(s) \leq \dist_{\Gamma}(s, u_0)\leq \diam(\Gamma)$. Since $s$ is arbitrary, we obtain $h_{LSG}(\Gamma) \leq \diam(\Gamma)$ by Definition~\ref{Def2.14New}\eqref{Def2.14-3}. Furthermore, if $\Gamma$ is connected, not LSG-complete, the case $h_{LSG}(\Gamma) \in \{0, 1\}$ is excluded by virtue of Item (ii). That means $2 \leq h_{LSG}(\Gamma), \, n \geq 3$. So, Assertion~\eqref{thm1-iii} holds. 
    \end{itemize}            
\end{itemize} Theorem~\ref{thm1} has been fully proven. 
\end{proof}

\subsection{A procedure for computing LSG-height}\label{Subsec2.3New}

Let $\Gamma$ be an $n$-LSG constructed on $\F^n = \s \{e_1, \ldots, e_n\}$. Due to Theorem~\ref{thm1}, its structure depends on the LSG-height $h_{LSG}(\Gamma)$. Take an arbitrary vertex $s \in \F^n_0$. Our objective is to compute exactly the LSG-neighborhood sequence $\Nn^k_{LSG}(s)$ for $k \geq 1$. For convenience, we denote $W_k \co \Nn^k_{LSG}(s)$.

\begin{enumerate}[\bf Step 1.]
    \item Compute $W_1 \co \Nn^1_{LSG}(s) = N_{LSG}(s)$. If $W_1 = \F s$ then $h_{LSG}(s) = 0$. If $\F s \neq W_1 = \F^n$ then $h_{LSG}(s) = 1$. Otherwise, go to Step 2 with $k = 1$. 

    \item Compute $W_{k+1}$ whenever $W_k = \Nn^k_{LSG}(s) = \s \{w_1, \dots, w_m\}$. By Definitions~\ref{Def2.9New} and~\ref{Def2.12New},
    \[
        W_{k+1} \co \Nn^{k+1}_{LSG}(s) \co  \Nn_{LSG}(\Nn^k_{LSG}(s)) = \sum_{w \in (W_k)_0}N_{LSG}(w).
    \]
    Since $S \co \sum \limits_{i = 1}^m N_{LSG}(w_i) \leq W_{k+1}$, one first computes $N_{LSG}(w)$ for an arbitrary element $w = t_1w_1 + \dots + t_mw_m \in (W_k)_0$. If there exist $(t_1, \dots, t_m)$ such that $N_{LSG}(w) \nleq S$, add the corresponding subspaces $N_{LSG}(w)$ to $S$. Otherwise, $W_{k+1} = S$.
    \begin{itemize}
         \item If $W_{k+1} = W_k$, then $h_{LSG}(s) = k$. 
         \item If $W_{k+1} = \F^n$, then $h_{LSG}(s) = k+1$.
         \item Otherwise replace $k$ by $k+1$ and restart this step.
    \end{itemize}
    By Proposition~\ref{Prop2.13New}, this step must stop since it is repeated at most $n-1$ times.
\end{enumerate}

\begin{example}\label{Ex2.17New} Consider $\F^3 = \s\{e_1, e_2, e_3\}$ and $\omega\colon \F^3 \times \F^3 \to \F$ as follows:
\[\omega\colon \F^3 \times \F^3 \to \F, \, (u, v) = \bigl((u_1, u_2, u_3), (v_1, v_2, v_3)\bigr) \mapsto \omega(u, v) \co u_1 v_2 - u_2 v_1.\]
Let $\Gamma = (V_{\Gamma}, E_{\Gamma})$ be the following graph
\begin{itemize}
    \item {\bf The vertex set} $V_{\Gamma} \co \F^3_0$.
    \item {\bf The edge set} $E_{\Gamma} \co \bigl{\{}\{u, v\}\mid u \ne v, \, \omega(u, v) = 0\bigr{\}}$.
\end{itemize}
For any $u,v \in \F^3$, $u \sim_{LSG} v$ if and only if $\omega(u,v)=0$. It is obvious that $\omega$ is bilinear and skew-symmetric, so by Definition~\ref{Def2.1New}, $\Gamma$ is a $3$-LSG constructed on $\F^3$.
    
To apply the above procedure, we take an arbitrary vertex $s = (s_1, s_2, s_3) \in \F^3_0$. In Step 1, we compute $W_1 = N_{LSG}(s) = \{u \in \F^3 \mid u \sim_{LSG} s\}$ which is exactly the solution space of the linear equation $\omega(s, u) = s_1u_2 - s_2u_1 = 0$. It always contains $s$ and $e_3$. There are two cases as follows.
    \begin{itemize}
        \item If $s_1 = s_2 = 0$: then $s = s_3e_3 \in \F_0 e_3$ 
        and the above equation holds automatically. Hence, $W_1 = N_{LSG}(s) = \F^3$ and $\F s \neq W_1$. So $h_{LSG}(s) = 1$.

        \item If $(s_1, s_2) \neq (0, 0)$ then $W_1 = N_{LSG}(s) = \s \{s, e_3\} \neq \F s$. So we turn to Step 2 with $W_1 = \s \{s, e_3\}$. Since $N_{LSG}(e_3) = \F^3$, we have
        \[
            \F^3 = N_{LSG}(s) + N_{LSG}(e_3) = S \leq W_2 \leq \F^3.
       \]
       This means $W_2 = \F^3$, and the second item of Step 2 implies $h_{LSG}(s) = 2$.
    \end{itemize}
    Since $s \in \F^3_0$ is arbitrary, we have $h_{LSG}(\Gamma) = \max \limits_{s \in \F^3_0} h_{LSG}(s) = 2$.  
\end{example}

\begin{example}\label{Ex2.18New} Let $\F = \R$.		Write every vector of $\R^4$ as $u=(a,b)$ with $a=(a_1,a_2)$, $b=(b_1,b_2) \in \R^2$. For $a \in \R^2$, put
	\[
	   D(a) \co a_1
       \begin{pmatrix}
        1 & 0 \\ 0 & 1
       \end{pmatrix}
       + a_2
       \begin{pmatrix}
        0 & -1 \\ 1 & 0
       \end{pmatrix}
       =
       \begin{pmatrix}
        a_1 & -a_2 \\ a_2 & a_1
       \end{pmatrix}.
	\]
For $u=(a,b), u'=(a',b') \in \R^4$, consider the map $\omega \colon \R^4 \times \R^4 \to \R^2$ as follows
	\begin{equation}\label{For2.8New}
		\omega(u,u')\co D(a)b'-D(a')b
			=\begin{pmatrix}a_1b'_1-a_2b'_2-a'_1b_1+a'_2b_2\\ a_2b'_1+a_1b'_2-a'_2b_1-a'_1b_2\end{pmatrix}
	\end{equation}
Let $\Gamma = (V_{\Gamma}, E_{\Gamma})$ be the following graph
\begin{itemize}
    \item {\bf The vertex set} $V_{\Gamma} \co \R^4_0$.
    \item {\bf The edge set} $E_{\Gamma} \co \bigl{\{}\{u, u'\}\mid u \ne u', \, \omega(u, u') = 0\bigr{\}}$.
\end{itemize}
In other words, two distinct vertices $u,u' \in \R^4_0$ are adjacent if and only if $\omega(u,u')=0$. It is obvious that $\omega$ is bilinear and skew-symmetric, so the conditions of Definition~\ref{Def2.1New} hold. Therefore, $\Gamma$ is a 4-LSG constructed on $\R^4$. 
		
To run the procedure above, we use the complex structure hidden in $D$. First, we consider the $\R$-linear isomorphism $j \colon \R^2 \to \C$ as follows
\[
    a = (a_1, a_2) \mapsto j(a) \co a_1 + i a_2, \qquad \text{($i$ is the imaginary unit)}.
\]
Then, for every $a = (a_1, a_2), \, b = (b_1, b_2) \in \R^2$, we have  $j\bigl(D(a)b\bigr)=j(a)\,j(b)$ and
\begin{equation}\label{For2.9New}
    j\bigl(\omega(u, u')\bigr) = j(a)j(b') - j(a')j(b); \quad \forall u = (a, b), u' = (a', b') \in \R^4.
\end{equation}
Next, we consider the $\R$-linear isomorphisms $\sigma \colon \R^4 \to \R^4$ and $\varphi\colon\R^4\to\C^2$ as follows
\[
    \begin{cases}
    (a, b) = (a_1, a_2, b_1, b_2) \mapsto \sigma (a_1, a_2, b_1, b_2) \co (-a_2, a_1, - b_2, b_1), \\
    (a, b) = (a_1, a_2, b_1, b_2) \mapsto \varphi(a, b) \co \bigl(j(a), j(b)\bigr).
\end{cases}
\]
It is clear that $\varphi\bigl(\sigma(s)\bigr)=i\,\varphi(s), \, \forall s \in \R^4$. Take an arbitrary vertex $s = (a_0, b_0) \in\R^4_0$, with $\varphi(s) = \bigl(j(a_0), j(b_0)\bigr) \ne (0, 0)$. There are two steps as follows. 
    \begin{itemize}
	\item \textbf{Step 1.} By equation~\eqref{For2.9New}, 
    \[W_1=N_{LSG}(s)=\bigl{\{} (a, b)\in\R^4 \equiv \C^2 \mid j(a_0)j(b) - j(a)j(b_0) = 0\bigr{\}}.\]   
    Since $\varphi(\sigma(s))=i\varphi(s)$ and $i\notin\R$, the vectors $s$ and $\sigma(s)$ are $\R$-linearly independent. By a direct calculation, we obtain
		\begin{equation}\label{For2.10New}
			W_1=N_{LSG}(s)=\varphi^{-1}\bigl(\C\varphi(s)\bigr)=\s\{s,\sigma(s)\}, \, \, \dim W_1=2 .
		\end{equation}
    Here, $\C u$ denotes the complex line through $u \in \C^2 \equiv \R^4$. Thus $\R s\neq W_1\neq\R^4$, and we turn to Step~2 with $k=1$.
	\item \textbf{Step 2.} Here, $W_1 = \s \{s, \sigma(s)\}$. 
    Let $w=t_1s+t_2\sigma(s)\in(W_1)_0$ be arbitrary, i.e., $(t_1,t_2)\neq(0,0)$. Then $\varphi(w)=\lambda\varphi(s)$ with $\lambda\co t_1+it_2\in\C_0$. For every $u = (a, b) \in \C^2 \equiv \R^4$, by Relation~\eqref{For2.9New}, we have
	\[
	j\bigl(\omega(w,u)\bigr)=\lambda \bigl[j(a_0) j(b) -j(a) j(b_0)\bigr] 
    =\lambda\, j\bigl(\omega(s,u)\bigr).
	\]
	Hence, $\omega(w,u)=0$ if and only if $\omega(s,u)=0$, i.e.,
		\begin{equation}\label{For2.11New}
			N_{LSG}(w)=N_{LSG}(s)=W_1; \, \, \forall w\in(W_1)_0.
		\end{equation}
	So $S\co N_{LSG}(s) + N_{LSG}(\sigma(s)) = W_1$, and there is no $w \in (W_1)_0$ such that $N_{LSG}(w) \nleq S$. Therefore, $W_2=\Nn^2_{LSG}(s)=W_1$, and the first item of Step~2 gives $h_{LSG}(s)=1$.
	\end{itemize}
	Since $h_{LSG}(s)=1$ for any $s\in\R^4_0$, we have $h_{LSG}(\Gamma) = \max \limits_{s \in \R^4_0} h_{LSG}(s) = 1$.
\end{example}

\section{The CG of a Lie algebra with trivial center viewed as an LSG}\label{Sec3New}

\subsection{Lie algebra graphs}\label{Subsec3.1New}

\begin{proposition and definition}\label{PD3.1New}
Let $\G$ be a Lie algebra with trivial center.
    \begin{enumerate}[(i)]
        \item\label{PD3.1-i} The CG $\Gamma(\G)$ is an LSG constructed on $\G$ with $\dim \Gamma(\G) = \dim (\G)$ which is called the Lie algebra graph (LAG) of $\G$. When $\dim(\G) = n \geq 2$, we say that $\Gamma(\G)$ is an $n$-LAG of $\G$.
        
        \item\label{PD3.1-ii} For every vertex $s$, $N_{LSG}(s) = C_{\G}(s) \co \{x \in \G \mid [s, x] = 0\}$, the centralizer of $s$.
    \end{enumerate} 
\end{proposition and definition}
\begin{proof}
Since $\Zz(\G)=\{\boldsymbol{0}\}$, $V(\G)= \G\setminus\Zz(\G) = \G\setminus\{\boldsymbol{0}\} \equiv \G_0$. Furthermore, 
\[E(\G) \co \bigl{\{}\{x, y\} \, |\, x, \, y \in V(\G), \, x \ne y, \, [x, y] = 0\bigr{\}}.\]
\begin{itemize}
    \item {\bf Proof of~\eqref{PD3.1-i}.} The bilinearity of the Lie bracket implies that $E(\G)$ is linearly stable, that is, $L_1$, $L_2$, and $L_3$ are simultaneously satisfied. In other words, Assertion~\eqref{PD3.1-i} holds.
    
    \item {\bf Proof of~\eqref{PD3.1-ii}.} Assertion~\eqref{PD3.1-ii} follows directly from the definition of the centralizer of every $s \in \G_0$, Definition~\ref{Def1.3New} and  Definition~\ref{Def2.2New}.
\end{itemize}
\end{proof}

\begin{remark}\label{Re3.2New} The following three observations are obvious yet important.
\begin{enumerate}[(1)]    
    \item\label{Re3.2-1} {\bf LAGs are compatible with LSG-direct sums}: Assume that $\G_1, \G_2$ are two non-zero Lie algebras of trivial centers. For every $u = (u_1, u_2), v = (v_1, v_2) \in \G_1 \oplus_{\text{Lie}}\G_2$, since $[u, v]_{\G_1 \oplus_{\text{Lie}}\G_2}=0$ if and only if $[u_1,v_1]_{\G_1}=0$ and $[u_2,v_2]_{\G_2}=0$, $\Gamma(\G_1\oplus_{\text{Lie}}\G_2)=\Gamma(\G_1)\oplus_{LSG}\Gamma(\G_2)$.
    
    \item\label{Re3.2-2} {\bf Each Lie monomorphism induces an LSG-monomorphism}: Let $\G, \G'$ be Lie algebras with trivial centers and $f \colon \G \hookrightarrow_{\text{Lie}} \G'$ be a Lie monomorphism. 
    Then for all $x, y \in \G$, $[x,y]_{\G}=0$ if and only if $[f(x),f(y)]_{\G'}=0$. By Definitions~\ref{Def2.4New} and~\ref{Def2.7New}, 
    the Lie monomorphism $f$ induces an LSG-monomorphism 
    $LSG(f) \colon \Gamma(\G) \hookrightarrow_{LSG} \Gamma(\G')$.

    \item\label{Re3.2-3} {\bf LSG-invariants of LAGs}. For any Lie algebra $\G$ with trivial center, the LSG-height and the LSG-colsed dimension of LAG $\Gamma(\G)$ are LSG-invariants. They are called {\it LSG-invariants of $\G$}.
\end{enumerate} 
\end{remark}

\begin{example}\label{Ex3.3New}
Below are two typical examples.
\begin{enumerate}[(1)]
    \item\label{Ex3.3-1} {\bf The real affine Lie algebra}. Let $\F = \R$.
    \[\af(\R)=\s\{e_1,e_2\}, \, \, [e_1,e_2]=e_2.\]
    The center of $\af(\R)$ is trivial and all centralizers of non-zero elements are $1$-dimensional, so $h_{LSG}\bigl(\Gamma(\af(\R))\bigr)=0$. By Theorem~\ref{thm1}(ii)(a), we have
    \[\Gamma(\af(\R)) \cong_{LSG} \bigsqcup_{\widetilde{s} \in (\s\{e_1, e_2\}\setminus\{0\})/\sim \, = \, \mathbb{P}^1_{\R}} \boldsymbol{K}_{\widetilde{s}}.\]
    
    \item\label{Ex3.3-2} {\bf The complex special linear Lie algebra of degree 2}. Let $\F = \C$.
    \[\mathfrak{sl}_2(\C)=\s\{x,y,z\}, \, \, [x,y]=z, \, \, [y,z]=2y, \, \, [x,z]=-2x.\]
    All centralizers of non-zero elements are $1$-dimensional, so $h_{LSG}\bigl(\Gamma(\mathfrak{sl}_2(\C))\bigr)=0$ and 
    \[\Gamma(\mathfrak{sl}_2(\C)) \cong_{LSG} \bigsqcup_{\widetilde{s} \in (\s\{x, y, z\}\setminus\{0\})/\sim \, = \, \mathbb{P}^2_{\C}}
    \boldsymbol{K}_{\widetilde{s}}.\]
\end{enumerate}
\end{example}

\subsection{LAGs of commutative transitive Lie algebras}\label{Subsec3.2New}
This entire subsection is devoted to the LAGs of a special class of Lie algebras known as CT ones. First, we introduce their definition.

\begin{definition}[{see \cite[Introduction]{Klep2010}}]\label{Def3.4New}
    Let $\G$ be an $n$-dimensional Lie algebra ($n \geq 1$). $\G$ is called  {\it commutative transitive} (or {\it CT}, for short) if either $\G$ is abelian, or $\G$ is non-abelian and for all $x, y, z \in \G_0$,
    \begin{equation}\label{For3.1New}
        \bigl([x, y] = 0\bigr) \wedge \bigl([y, z] = 0\bigr) \Longrightarrow \bigl([x, z] = 0\bigr)
    \end{equation}                       
\end{definition}

The class of CT Lie algebras (CTLA) was first addressed by I. Klep and P. Moravec~\cite[Introduction]{Klep2010} in 2010. By~\cite[Lemma 3]{Gorb2017}, every non-abelian CTLA $\G$ has trivial center: $\Zz(\G) = \{\boldsymbol{0}\}$. Consequently, the CG $\Gamma(\G)$ of $\G$ is an LAG. Furthermore, the class of all non-abelian CTLAs over an algebraically closed field of characteristic $0$ was classified by I. Klep and P. Moravec~\cite{Klep2010} in 2010. Specifically, we have the following result.

\begin{proposition}[{\cite[Theorems 3, 6, 7]{Klep2010}}]\label{Prop3.5New}
Let $\F$ be algebraically closed of characteristic $0$ and $\G$ a non-abelian CTLA of dimension $n\ge2$. Then $\G$ is simple or solvable. Furthermore, we have
\begin{enumerate}[(i)]
    \item If $\G$ is simple then $\G \cong \mathfrak{sl}_2(\F)$.
    
    \item If $\G$ is solvable then $\G \cong \F^{n-1} \oplus_D \F e_0$ with $D=\ad_{e_0}|_{\F^{n-1}} \in GL_{n-1}(\F)$. \hfill $\square$
\end{enumerate}
\end{proposition}

\begin{theorem}[{{\bf LSG-height and LAGs of CTLAs}}]\label{thm2}
Let $\G$ be an $n$-dimensional non-abelian Lie algebra with trivial center ($n \geq 2$) and $\Gamma(\G)$ be its $n$-LAG. The following are equivalent:
\begin{enumerate}[(i)]
    \item\label{thm2-i} $h_{LSG}\bigl(\Gamma(\G)\bigr) \leq 1$.
    
    \item\label{thm2-ii} $C_{\G}(s)$ is an abelian subalgebra of $\G$ for every $s \in \G_0$.
    
    \item\label{thm2-iii} $\G$ is a CTLA.
\end{enumerate}
In that case, $\s_{LSG}(s) = C_{\G}(s)$ for every $s \in \G_0$, the LSsubG-connected components of $\Gamma(\G)$ are exactly the LSG-completes $\boldsymbol{K}_{C_{\G}(s)_0}$, the $C_{\G}(s)$ are exactly the maximal abelian subalgebras of $\G$ and $\Gamma(\G)$ is always disconnected. Moreover, 
\begin{equation}\label{For3.2New}
    (h_{LSG}\bigl(\Gamma(\G)\bigr) = 0) \Longleftrightarrow \bigl[(C_{\G}(s) = \F s; \, \forall s \in \G_0) \wedge \bigl(\Gamma(\G)\cong_{LSG}\bigsqcup_{\widetilde{s} \in \mathbb{P}^{n-1}_{\F}} \boldsymbol{K}_{\widetilde{s}}\bigr)\bigr]
\end{equation}
\end{theorem}

\begin{proof}[{\bf Proof of Theorem~\ref{thm2}}]
By Proposition and Definition~\ref{PD3.1New}(ii) and the Jacobi identity, we have
\begin{equation}\label{For3.3New}
    N_{LSG}(s)=C_\G(s) = \Nn^1_{LSG}(s) \leq \Nn^2_{LSG}(s)=\sum_{w\in C_\G(s)_0}C_\G(w); \forall s \in \G_0.
\end{equation}
So for every $s \in \G_0$, we also have
\begin{equation}\label{For3.4New}
    \bigl(h_{LSG}(s)\le1\bigr) \iff \bigl(C_\G(w)\subseteq C_\G(s); \forall w\in C_\G(s)_0\bigr).
\end{equation}
\begin{itemize}
    \item {\bf Proof of~\eqref{thm2-i} $\Rightarrow$~\eqref{thm2-ii}}:
    Assume that $h_{LSG}\bigl(\Gamma(\G)\bigr) \leq 1$, so that $h_{LSG}(s) \leq 1$ for every $s \in \G_0$. By Relations~\eqref{For3.3New} and ~\eqref{For3.4New}, we have
    \begin{equation}\label{For3.5New}
    C_{\G}(w) \subseteq C_{\G}(s); \, \, \forall w \in C_{\G}(s)_0, \, \, \forall s \in \G_0
    \end{equation}
    Now fix $s \in \G_0$ and let $w \in C_{\G}(s)_0$. Then $s \in C_{\G}(w)_0$, so Relation~\eqref{For3.5New} applied to the vertex $w$ (which is allowed, since Relation~\eqref{For3.4New} holds for every vertex) gives $C_{\G}(s) \subseteq C_{\G}(w)$. That means 
    \begin{equation}\label{For3.6New}
    C_{\G}(w) = C_{\G}(s); \, \, \forall w \in C_{\G}(s)_0 .
    \end{equation}
    Let $w, w' \in C_{\G}(s)_0$. Of course,  
    $w' \in C_{\G}(s) = C_{\G}(w)$ by Formula~\eqref{For3.6New}, i.e.\ $[w, w'] = 0$. Hence $C_{\G}(s)$ is abelian, and (ii) holds.
    
    \item {\bf Proof of~\eqref{thm2-ii} $\Rightarrow$~\eqref{thm2-i}}:   Assume that $C_{\G}(s)$ is abelian for every $s \in \G_0$. Fix $s \in \G_0$ and let $w \in C_{\G}(s)_0$. Since $C_{\G}(s)$ is an abelian subalgebra containing $w$, every element of $C_{\G}(s)$ commutes with $w$, i.e.\ $C_{\G}(s) \subseteq C_{\G}(w)$. Symmetrically, $C_{\G}(w)$ is abelian and contains $s$, so $C_{\G}(w) \subseteq C_{\G}(s)$. Hence Relation~\eqref{For3.6New} holds, and by Relation~\eqref{For3.4New}, we get
    \[
    \Nn^2_{LSG}(s) = \sum_{w \in C_{\G}(s)_0} C_{\G}(w) = C_{\G}(s) = \Nn^1_{LSG}(s).
    \]
    Thus the LSG-neighborhood sequence of $s$ is stationary from the first step on, i.e.\ $h_{LSG}(s) \leq 1$. Since $s \in \G_0$ was arbitrary, $h_{LSG}\bigl(\Gamma(\G)\bigr) \leq 1$ and (i) holds.
    
    \item {\bf Proof of~\eqref{thm2-ii} $\Rightarrow$~\eqref{thm2-iii}}:
    Assume that $C_{\G}(s)$ is abelian for every $s \in \G_0$. Let $x, y, z \in \G$ with $y \neq \boldsymbol{0}$, $[x, y] = 0$ and $[y, z] = 0$. Then $x \in C_{\G}(y)$ and $z \in C_{\G}(y)$, because the bracket is alternating. Since $C_{\G}(y)$ is abelian, $[x, z] = 0$. Hence $\G$ is CT.
    
    \item {\bf Proof of~\eqref{thm2-iii} $\Rightarrow$~\eqref{thm2-ii}}:
    Assume that $\G$ is CT, that is, for all $x, y, z \in \G$ with $y \neq 0$,
    \begin{equation}\label{For3.7New}
    \bigl([x, y] = 0 \ \text{and} \ [y, z] = 0\bigr) \ \Longrightarrow \ [x, z] = 0 .
    \end{equation}
    Fix $s \in \G_0$. We first note that $C_{\G}(s)$ is a subalgebra of $\G$: if $x, z \in C_{\G}(s)$, then by the Jacobi identity $[s, [x, z]] = [[s, x], z] + [x, [s, z]] = 0$, so $[x, z] \in C_{\G}(s)$.
    \item[] To show that $C_{\G}(s)$ is abelian, let $x, z \in C_{\G}(s)$ be arbitrary. We apply Relation~\eqref{For3.7New} to the triple $(x, y, z)$ with $y \co s$, and check its three hypotheses.
    \begin{itemize}
        \item $y = s \neq \boldsymbol{0}$, since $s \in \G_0$.
        \item $[x, y] = [x, s] = 0$, since $x \in C_{\G}(s)$.
        \item $[y, z] = [s, z] = -[z, s] = 0$, since $z \in C_{\G}(s)$ and the bracket is alternating.
    \end{itemize}
    Hence Relation~\eqref{For3.7New} yields $[x, z] = 0$. As $x, z \in C_{\G}(s)$ were arbitrary, any two elements of $C_{\G}(s)$ commute, i.e.\ $C_{\G}(s)$ is an abelian subalgebra of $\G$. Since $s \in \G_0$ was arbitrary, (ii) holds.
    
    \item {\bf Proof of the additional assertions}: Suppose one of the statements (i), (ii), or (iii) is true; it follows that all three are true. For $s \in \G_0$, the LSG-neighborhood sequence in Definition~\ref{Def2.12New} is stationary from the first step on, so $\s_{LSG}(s) = \Nn^1_{LSG}(s) = C_{\G}(s)$ by Definition~\ref{Def2.14New}. By Theorem~\ref{thm1}, Items~(ii)(a) and~(ii)(b), every connected component of $\Gamma(\G)$ is the complete LSG
    $\boldsymbol{K}_{\s_{LSG}(s)\setminus\{\boldsymbol{0}\}} = \boldsymbol{K}_{C_{\G}(s)_0}$ for some $s \in \G_0$, and $\Gamma(\G)$ is connected only if it is the complete LSG $\boldsymbol{K}_{\G_0}$; the latter would mean
    $[\G, \G] = \{\boldsymbol{0}\}$, which is excluded, so $\Gamma(\G)$ is disconnected. Finally, an abelian subalgebra containing $s$ is contained in $C_{\G}(s)$, which is abelian by (ii); hence $C_{\G}(s)$ is the unique maximal abelian subalgebra containing $s$, and every maximal abelian subalgebra $A$ equals $C_{\G}(s)$ for any $s \in A_0$. The statement on LSG-height $0$ is the case $C_{\G}(s) = \Nn^1_{LSG}(s) = \F s$ for all $s$, and by Theorem~\ref{thm1}\eqref{thm1-ii}(a), Relation~\eqref{For3.2New} is correct.
\end{itemize}
This completes the proof of Theorem~\ref{thm2}.
\end{proof}

\begin{example}[{\bf CT Lie algebras of LSG-height $1$}]\label{Ex3.6New} Regarding non-abelian CT Lie algebras with LSG-height $1$, we have the two typical examples below.
\begin{enumerate}[(1)]
    \item\label{Pro3.6New-1} {\bf A $3$-dimensional complex almost abelian algebra (Kolman~\cite{Kol1965})}. Let $\F = \C$. We consider the following $3$-dimensional complex Lie algebra
    \[
    \G_1=\C^2\oplus_{\id}\C e_0 \co \s\Bigl{\{}e_1, e_2; e_0 \, |\, [e_0,e_1]=e_1, \, \, [e_0,e_2]=e_2\Bigr{\}}.
    \]     
    For any $x=u+\alpha e_0, \, y=v+\beta e_0 \in \G_1$ with $u, v \in \s\{e_1, e_2\} \equiv \C^2$, we have $[x,y]=\alpha v-\beta u$. Hence $C_{\G_1}(x)=\C^2$ if $x\in\C^2_0$, and $C_{\G_1}(x)=\C x$ if $x\notin\C^2$. All centralizers are abelian, so $\G_1$ is CT by Theorem~\ref{thm2}, with $h_{LSG}(\Gamma(\G_1))=1$ and $\dim_{LSG}(\Gamma(\G_1))=2=\dim\Gamma(\G_1)-1$. The LAG of $\G_1$ is given by Relation~\eqref{For2.5New} in Item~\eqref{thm1-ii}(b) of Theorem~\ref{thm1} applied for $n = 3$. Namely,
    \[
    \Gamma(\G_1) \cong_{LSG} \boldsymbol{K}_{\C^2_0} \sqcup \Bigl(\bigsqcup_{\widetilde{s}\in \bigl(\C^3_0\setminus\C^2_0\bigr)/\sim \,=\,\mathbb{P}^2_{\C}\setminus \mathbb{P}^1_{\C}}\boldsymbol{K}_{\widetilde{s}} \Bigr)
    \]  
    
    \item\label{Pro3.6New-2} {\bf A $4$-dimensional real CTLA which is not almost abelian.} Let $\F = \R$, $n = 4$. We consider the following $4$-dimensional real Lie algebra 
    \[ \G_2= \s\Bigl{\{}x_1, x_2, y_1, y_2\, |\, [x_1,y_1]=y_1, \, [x_1,y_2]=y_2, \, [x_2,y_1]=y_2, \, [x_2,y_2]=-y_1\Bigr{\}}. 
    \] 
    That means $\G_2=\R^2_y\oplus_D\R^2_x = \{(a, b)\, |\, a = (a_1, a_2) \in \R^2_x, \, b = (b_1, b_2) \in \R^2_y\}$ with $\R^2_x \co \s\{x_1, x_2\}$, $\R^2_y \co \s\{y_1, y_2\}$; the linear map $D$ is given by
    \[ D \colon \R^2_x \to \mathfrak{gl}_2(\R) \equiv \Der(\R^2_y); \, a = (a_1, a_2) \mapsto D(a) \co a_1 \begin{pmatrix}1&0\\0&1\end{pmatrix} + a_2 \begin{pmatrix}0&-1\\1&0\end{pmatrix}. 
    \] 
    For any $u = (a, b) = (a_1, a_2, b_1, b_2), u' = (a', b') = (a'_1, a'_2, b'_1, b'_2) \in \G_2$, we get
    \[[u, u'] \co D(a)(b') - D(a') (b) = \begin{pmatrix} a_1b'_1 - a_2b'_2 - a'_1b_1 + a'_2b_2\\ a_1b'_2 + a_2b'_1 - a'_1b_2 - a'_2b_1 \end{pmatrix} \in \R^2_y.
    \] 
    This is, in fact, the $4$-dimensional real Lie algebra whose LAG is exactly the LSG considered in Example~\ref{Ex2.18New}.
    By a direct calculation, we obtain $[\G_2, \G_2] = \R^2_y$, the center $\Zz(\G_2) = \{\boldsymbol{0}\}$ and $\G_2$ is a $4$-dimensional CTLA. In particular, $\Gamma(\G_2)$ is the $4$-LAG of height $h_{LSG}\bigl(\Gamma(\G_2)\bigr) \leq 1$. Moreover, for any $u = (a, b)\in (\G_2)_0$, we have  
    \[ \R u < C_{\G_2}(u) =
    \begin{cases}
        \Bigl{\{} u' = (a', b') \in \G_2 \,|\, D(a) b' = D(a')(b)\Bigr{\}}\cong \R^2 \quad \text{if} \, \, \, a \ne 0,\\
        \R^2_y \quad \text{if} \, \, \, a = 0 \ne b.
    \end{cases}  
    \] 
    By Theorem~\ref{thm2} (and as determined in Example~\ref{Ex2.18New}), we get $h_{LSG}\bigl(\Gamma(\G_2)\bigr) = 1, \, \dim_{LSG}\bigl(\Gamma(\G_2)\bigr) = 2 = \dim\Gamma(\G_2) - 2$. If $A$ is an abelian subalgebra of $\G_2$ and $\boldsymbol{0}\ne u\in A$, then $A \subseteq C_{\G_2}(u)$. Therefore, $\dim(A) \leq 2$. Hence $\G_2$ has no abelian ideal of codimension $1$, i.e. $\G_2$ is not almost abelian (see Definition~\ref{Def4.7New}). Furthermore, $\G_2$ is a CTLA of LSG-height $1$ with the LSG-closed dimension $\dim_{LSG}\bigl(\Gamma(\G_2)\bigr)$ strictly between $1$ and $\dim(\G_2) - 1$. In other words, Case~3 of Theorem~\ref{thm1}\eqref{thm1-ii}(b) does occur for LAGs over $\R$. It should also be noted that $\G_2$ is a $4$-dimensional real CTLA and, naturally, does not appear in the classification of I. Klep and P. Moravec~\cite{Klep2010} in 2010 (see Proposition~\ref{Prop3.5New}), since $\R$ is not an algebraically closed field.
    \vskip 0.2cm
    \item[] Now, we will describe the LAG $\Gamma(\G_2)$ explicitly.  
    By Theorem~\ref{thm2} and Definition~\ref{Def2.4New}, all LSsubG-connected components of $\Gamma(\G_2)$ are the LSG-completes $\boldsymbol{K}_{(C_{\G_2}(s))_0}$ of dimension $2$ with $s \in (\G_2)_0$. Furthermore, $\Gamma(\G_2)$ has infinitely many LSsubG-connected components of dimension $\dim_{LSG}(\Gamma(\G_2)) = 2$. Of course, it does not contradict the uniqueness statement in the third case of Theorem~\ref{thm1}\,(ii)(b), since here $2\dim_{LSG}(\Gamma(\G_2)) = 4 = n$, so the assumption $2\dim_{LSG}(\Gamma) > n$ does not hold. Now, we need to find a way to ``count'' the connected components of $\Gamma(\G_2)$. 

    First, we identify $\R^2$ with $\C$ (viewed as real vector spaces) via a linear isomorphism $j \colon \R^2 \to \C, \, a = (a_1, a_2) \mapsto j(a) \co a_1 + i a_2$, where $i$ is the imaginary unit. In that case, $\G_2 \equiv \C^2$ (as $4$-dimensional real vector spaces) and we obtain 
    \begin{alignat*}{3}
        D\colon \R^2_x \equiv \C \to & \,  \mathfrak{gl}_2(\R) \subseteq \Der(\R^2_y) \equiv \text{End}_{\R}(\C)\\
        a = (a_1, a_2) &\mapsto D(a)(.) = a_1 \begin{pmatrix}
            1&0 \\ 0&1
        \end{pmatrix} + a_2 
        \begin{pmatrix}
            0&-1 \\ 1&0
        \end{pmatrix} = (a_1 + i a_2)(.)\\
        D(a) & \colon \C \to \C, \, \, b \mapsto D(a)(b) \co a\,b
    \end{alignat*}
    It means that, for every $u = (a, b), u' = (a', b') \in (\G_2)_0$, we have
    \begin{alignat*}{2}
        \bigl(u' \in N_{LSG}(u) = C_{\G_2}(u)\bigr) &\Leftrightarrow \bigl([u, u'] = 0\bigr)\\
        &\Leftrightarrow \, (a\,b' = a'\,b) \Leftrightarrow \bigl[(a:b) = (a':b') \in \C^2_0/\sim_{\C} \bigr].
    \end{alignat*}    
    Here $\sim_{\C}$ denotes the proportionality relation in $\G_2 \equiv \C^2$. In other words, $C_{\G_2}(s) = \C s$ for any $s \in (\G_2)_0$ and the family of LSsubG-connected components of $\Gamma(\G_2)$ is indexed by $(\G_2)_0/\sim_{\C} = \mathbb{P}^1_{\C}$. For any $s \in (\G_2)_0$, as is typical, the notation $\widetilde{s} \co \C_0 s$ indicates the equivalence class of $s$. Thus, the LAG $\Gamma(\G_2)$ is described as follows:  
    \[ \Gamma(\G_2) \cong_{LSG} \bigsqcup_{\widetilde{s} \in \mathbb{P}^1_{\C}} \boldsymbol{K}_{\widetilde{s}}. 
    \]

\end{enumerate}
\end{example}

\section{The reduced commuting graph of a Lie algebra}\label{Sec4New}

\subsection{Definition, well definedness and invariance under Lie isomorphisms}\label{Subsec4.1New}

\begin{definition}[{{\bf Reduced commuting graph - RCG}}]\label{Def4.1New}
Let $\G$ be a non-abelian Lie algebra with the center $\Zz(\G)$ ($\G \ne \Zz(\G) \triangleleft \G$). We fix a complementary subspace $\G_R$ of the center, i.e.\ $\G=\G_R\oplus\Zz(\G)$ and $\G_R \leq \G$ as a vector subspace. Then the subgraph $\Gamma_R(\G)=(V_R(\G),E_R(\G))$ of $\Gamma(\G)$ induced by $(\G_R)_0$ has 
\begin{itemize}
    \item {\bf The vertex set}: $V_R(\G)\co (\G_{R})_0$.
    \item {\bf The edge set}: $E_R(\G)\co\bigl{\{}\{x,y\}\mid x\neq y\in V_R(\G),\ [x,y]=0\bigr{\}}$.
\end{itemize}
Due to the bilinearity of the Lie bracket, $\Gamma_R(\G)$ is naturally an LSG constructed on $\G_R$ and it is called the {\it reduced commuting graph} (or {\it RCG}, for short) of $\G$. If $\G$ is $n+m$-dimensional and the center $\Zz(\G)$ is $m$-dimensional ($n, m \in \N; \, n \geq 1$), then $\G_R$ is always $n$-dimensional and $\Gamma_R(\G)$ is an $n$-LSG. So we call the RCG $\Gamma_R(\G)$ an {\it $n$-RCG} of $\G$.
\end{definition}

\begin{remark}\label{Re4.2New}
\begin{enumerate}[(1)]
    \item\label{Re4.2New-1} Since $\Gamma_R(\G)$ is an LSG, all notions of Section~\ref{Sec2New} (LSG-neighborhood, LSG-height, LSG-closed subspaces, LSG-isomorphisms, ...) can apply to it.
    
    \item\label{Re4.2New-2} Formally $\Gamma_R(\G)$ depends on $\G_R$. However, Theorem~\ref{thm3} below will show that it is uniquely determined up to an LSG-isomorphism.
    
    \item\label{Re4.2New-3} If $\Zz(\G)=0$, then $\G_R=\G$ and $\Gamma_R(\G)=\Gamma(\G)$ is an LAG.
    
    \item\label{Re4.2-4} If $\dim\Zz(\G)=m$, $\dim\G=m+n$ ($m, n \in \N, n \geq 1$) and $(x_1, \cdots, x_n, z_1, \cdots, z_m)$ is a basis of $\G$ with $\Zz(\G) = \s\{z_1, \cdots, z_m\}$, then a practical choice of $\G_R$ is $\G_R\co\s\{x_1,\dots,x_n\}$.
\end{enumerate}
\end{remark}

\begin{theorem}\label{thm3}
The RCG is both well defined and an invariant under Lie isomorphisms. Specifically, we have the following assertions. 
\begin{enumerate}[(i)]
    \item\label{thm3-i} {\bf Well definedness of the RCG}. If $\G$ is a non-abelian Lie algebra and $\G_{R},\G_{R'}$ are two complements of $\Zz(\G)$ in $\G$, then $\Gamma_{R}(\G)\cong_{LSG}\Gamma_{R'}(\G)$. In other words, $\Gamma_{R}(\G)$ is uniquely determined up to an LSG-isomorphism.
    \item\label{thm3-ii} {\bf The RCG is an invariant under Lie isomorphisms}. Let $f\colon\G\overset{\cong}{\to}_{Lie}\G'$ be a Lie isomorphism such that $f(\G_R)=\G'_R$. Then, $f$ always induces an LSG-isomorphism $\Gamma_R(f)\colon\Gamma_R(\G)\overset{\cong}{\to}_{LSG}\Gamma_R(\G')$. Hence, up to an LSG-isomorphism (whatever complements are chosen), the RCG is an {\it invariant} of the Lie isomorphism class of $\G$.
\end{enumerate}	
\end{theorem}

\begin{proof}[{\bf Proof of Theorem~\ref{thm3}}] We will prove each assertion.
\begin{itemize}
    \item {\bf Proof of~\eqref{thm3-i}.} Let $\pr'\colon\G=\G_{R'}\oplus\Zz(\G)\twoheadrightarrow \G_{R'}$ be the natural projection and $\varphi\co\pr'|_{\G_R}$. Since $\ker\pr'=\Zz(\G)$ and $\G_R\cap\Zz(\G)=\{\boldsymbol{0}\}$, $\varphi$ is injective. $\varphi$ is also surjective because $\dim\G_R=\dim\G_{R'}$. Thus, $\varphi\colon\G_R\overset{\cong}{\to}\G_{R'}$ is a linear isomorphism. For $u,v\in\G_R$ we have $u=\varphi(u)+z_u$, $v=\varphi(v)+z_v$ with $z_u,z_v\in\Zz(\G)$, whence $[u,v]=[\varphi(u),\varphi(v)]$. Thus $\{u,v\}\in E_R(\G)\iff\{\varphi(u),\varphi(v)\}\in E_{R'}(\G)$, and $\varphi$ induces an LSG-isomorphism.
    
    \item {\bf Proof of~\eqref{thm3-ii}.} Since $f\colon\G\overset{\cong}{\to}_{Lie}\G'$ is a Lie isomorphism with $f(\G_R)=\G'_R$, we obtain a linear isomorphism $f_R \co f\mid_{\G_R} \colon \G_R \overset{\cong}{\to}\G'_R$ which induces a bijection $f_{R_0}\colon V_R(\G) \to V_R(\G')$. For any $u, v \in V_R(\G)$ and $u \ne v$, we have 
    \begin{alignat*}{2}
        \{u, v\} \in E_R(\G) &\Leftrightarrow u\sim_{LSG} v \Leftrightarrow [u, v] = 0 \Leftrightarrow [f(u), f(v)] = 0\\
        & \Leftrightarrow [f_{R_0}(u), f_{R_0}(v)] = 0\\
        & \Leftrightarrow f_{R_0}(u) \sim_{LSG} f_{R_0}(v) \Leftrightarrow \{f_{R_0}(u), f_{R_0}(v)\} \in E_R(\G').
    \end{alignat*}
    Thus, $f_R \colon \G_R \overset{\cong}{\to}\G'_R$ defines an LSG-isomorphism
    \[\Gamma_R(f) \co f_{R_0}\colon\Gamma_R(\G)\overset{\cong}{\to}_{LSG}\Gamma_R(\G').\] 

    \item[] Finally, let $f\colon\G\overset{\cong}{\to}_{Lie}\G'$ be an arbitrary Lie isomorphism which does not necessarily satisfy the condition $f(\G_R) = \G'_R$ and $\G'_R$ an arbitrary complement of $\Zz(\G')$ in $\G'$. Since $f(\Zz(\G)) = \Zz(\G')$, the subspace $f(\G_R)$ is also a complement of $\Zz(\G')$ in $\G'$. Applying what has just been proven to $f$ and the complement $f(\G_R)$, and then~\eqref{thm3-i} to the two complements $f(\G_R)$ and $\G'_R$ of $\Zz(\G')$, we obtain $\Gamma_R(\G)\cong_{LSG}\Gamma_{f(\G_R)}(\G')\cong_{LSG}\Gamma_R(\G')$, where $\Gamma_{f(\G_R)}(\G')$ denotes the RCG of $\G'$ constructed on the complement $f(\G_R)$ of $\Zz(\G')$. Hence, up to an LSG-isomorphism (whatever complements are chosen), the RCG is an {\it invariant} of the Lie isomorphism class of $\G$.
\end{itemize} Theorem~\ref{thm3} is fully proven.
\end{proof}

Note that the converse of Theorem~\ref{thm3}\eqref{thm3-ii} does not hold. The example below shows that two non-isomorphic Lie algebras (which even have different dimensions) can still share the same RCG.

\begin{example}[{\bf Two non-isomorphic Lie algebras with LSG-isomorphic RCGs}]\label{Ex4.3New}
Let $\F = \R$.
\begin{enumerate}[(1)]
    \item\label{Ex4.3-1} {\bf The ordinary $3$-dimensional real space with the vector product}. Let
    \[
        \mathfrak{so}(3) \co \s\Bigl{\{}e_1, e_2, e_3 \,\mid\, [e_1, e_2] = e_3, \, [e_2, e_3] = e_1, \, [e_3, e_1] = e_2\Bigr{\}},
    \]
    i.e.\ $\R^3$ with the vector product. Since $[u, v] \co u \times v = 0$ if and only if $u, v$ are proportional, we have $\Zz(\mathfrak{so}(3)) = \{\boldsymbol{0}\}$. In particular, the LAG $\Gamma(\mathfrak{so}(3))$ is also the RCG $\Gamma_R(\mathfrak{so}(3))$ by Remark~\ref{Re4.2New} (3). 
    Furthermore, $N_{LSG}(v) = \R v$ for any $v \in \mathfrak{so}(3)\setminus\{\boldsymbol{0}\}$, and hence $h_{LSG}\bigl(\Gamma(\mathfrak{so}(3))\bigr) = 0$. As a direct consequence of Theorem~\ref{thm1}(ii)(a), we get
    \[
        \Gamma_R(\mathfrak{so}(3)) = \Gamma(\mathfrak{so}(3)) \cong_{LSG} \bigsqcup_{\widetilde{s} \in (\s\{e_1,e_2,e_3\}\setminus\{\boldsymbol{0}\})/\sim = \mathbb{P}^2_{\R}} \boldsymbol{K}_{\widetilde{s}}.
    \]
    
    \item\label{Ex4.3-2} {\bf The real $4$-dimensional diamond algebra}. Let
    \[
        \D \co \s\bigl\{x, y, z, t \,\big|\, [x, y] = z, \, [t, x] = -x, \, [t, y] = y\bigr\}.
    \]
    Then $\Zz(\D) = \R z$ and $\D_R \co \s\{x, y, t\}$ is a complement of the center. A direct computation shows that two elements of $\D_R$ commute if and only if they are proportional, so again $N_{LSG}(s) = \R s$ for every $s \in V_R(\D)$, $h_{LSG}(\Gamma_R(\D)) = 0$ and
    \[
        \Gamma_R(\D) \cong_{LSG} \bigsqcup_{\widetilde{s} \in (\s\{x,y,t\}\setminus\{\boldsymbol{0}\})/\sim = \mathbb{P}^2_{\R}} \boldsymbol{K}_{\widetilde{s}}.
    \]
\end{enumerate}
Consequently, $\Gamma_R(\mathfrak{so}(3)) \cong_{LSG} \Gamma_R(\D)$, although $\mathfrak{so}(3)$ and $\D$ are not Lie isomorphic because they have different dimensions.
\end{example}

\begin{remark}\label{Re4.4New} Naturally, we can make the following observations.
\begin{enumerate}[(1)]    
    \item\label{Re4.4-1} {\bf Structural advantage of RCGs}. For any Lie algebra $\G$, while the CG $\Gamma(\G)$ of $\G$ can be unwieldy due to the ``thickness'' of the center, the RCG is not only advantageous because of its linear structure, but also provides a ``streamlined''
    version that preserves all essential commutative information (see Theorem~\ref{thm4} in the next subsection) because it ``lies'' in
    a lower-dimensional space.

    \item\label{Re4.4-2} {\bf LSG-invariants of RCGs}. Of course, for any Lie algebra $\G$, the LSG-height and the LSG-closed dimension of RCG $\Gamma_R(\G)$ are LSG-invariants. They are called {\it LSG-invariants of $\G$}.
    
    \item {\bf Functorial property of RCGs}. The correspondence between Lie algebras and their RCGs can be described in the language of category theory. 
    \begin{itemize}
        \item\label{Re4.4-3} {\bf The category $\boldsymbol{\mathcal{R}\mathcal{L}_{iso}}$ of Lie algebras and reduced complements}  
        \begin{itemize}
            \item[$*$] {\bf Ob($\boldsymbol{\mathcal{R}\mathcal{L}_{iso}}$)} is the class of pairs $(\G, \G_R)$ where $\G$ is a non-abelian Lie algebra and $\G_R$ is a vector subspace in $\G$ such that $\G_R \oplus \Zz(\G) = \G$.
            \item[$*$] {\bf Mor}${}_{\boldsymbol{\mathcal{R}\mathcal{L}_{iso}}}$$\bigl((\G, \G_R), (\G', \G' _R)\bigr)$ is the set of Lie isomorphisms $f \colon \G \overset{\cong}{\to}_{Lie}\G'$ such that $f_R \coloneqq f|_{\G_R} \colon \G_R \overset{\cong}{\to} \G' _R$ is a linear isomorphism.
        \end{itemize}
        
        \item {\bf The category $\boldsymbol{\mathcal{LSG}_{iso}}$ of LSGs}
        \begin{itemize}
            \item[$*$] {\bf Ob($\boldsymbol{\mathcal{LSG}_{iso}}$)} is the class of LSGs $\Gamma$ constructed on any vector space $V$.
            \item[$*$] {\bf Mor}${}_{\boldsymbol{\mathcal{LSG}_{iso}}}$$(\Gamma, \Gamma')$ is the set of LSG-isomorphisms $LSG(f) \colon \Gamma  \overset{\cong}{\to}_{LSG}\Gamma'$ induced by linear isomorphisms $f \colon V \overset{\cong}{\to} V'$, where $\Gamma, \Gamma'$ are the LSGs constructed on the vector spaces $V, V'$, respectively.
        \end{itemize}
        \item[] According to Definition~\ref{Def4.1New}, we have the following functor: 
        \[\Gamma_R \colon \boldsymbol{\mathcal{R}\mathcal{L}_{iso}} \to \boldsymbol{\mathcal{LSG}_{iso}}; \quad (\G, \G_R) \mapsto \Gamma_R(\G), \, f \mapsto \Gamma_R(f) \co f_{R_0}.\]
        \item[] Furthermore, this functor is compatible with direct sums, meaning it maps the direct sum in $\boldsymbol{\mathcal{R}\mathcal{L}_{iso}}$ to the LSG-direct sum in $\boldsymbol{\mathcal{LSG}_{iso}}$. 
    \end{itemize}
\end{enumerate}    
\end{remark}

\subsection{RCG is a faithful replacement for CG}\label{Subsec4.2New}

Fix a Lie algebra $\G$. The CG $\Gamma(\G) = \bigl(V(\G), E(\G)\bigr)$ provides a graphical approach to the (generally non-simple) Lie structure of $\G$, while the RCG $\Gamma_R(\G) = \bigl(V_R(\G), E_R(\G)\bigr)$ has the advantage that its vertex set carries a natural linear structure, which interacts well with the edge set. This raises a fundamental question: Is the RCG $\Gamma_R(\G)$ a useful version of the CG $\Gamma(\G)$, i.e., do connectivity, the number of connected components and the diameter coincide for both graphs? As mentioned in the Introduction, the answer is affirmative (Theorem~\ref{thm4}); the key tool is the following result.

\begin{proposition}\label{Prop4.5New}

Let $\G$ be a (non-abelian) Lie algebra with the center $\Zz(\G)$ and the CG $\Gamma(\G)$. Assume that $\G_R$ is a complementary subspace of $\Zz(\G)$ in $\G$ and $\Gamma_R(\G)$ is the RCG of $\G$. The following statements hold.
	\begin{enumerate}[(i)]
		\item\label{Prop4.5-i} For every vertex $a \in V(\G)$, there exist unique $a_r \in V_R(\G)$ and $a_z \in \Zz(\G)$ such that $a = a_r + a_z$. Furthermore, for any distinct vertices $a, b \in V(\G)$, where $a = a_r + a_z, \, b = b_r + b_z$ with $a_r, b_r \in \G_R$, $a_z, b_z \in \Zz(\G)$, we always have 
        \begin{equation}\label{For4.1New}
        [a, b] = [a_r, b_r]; \quad \bigl(\{a,b\}\in E(\G)\bigr) \Leftrightarrow \Bigl[\bigl(a_r=b_r\bigr)\vee\bigl(\{a_r,b_r\}\in E_R(\G)\bigr)\Bigr]
        \end{equation}
        \item[] The map $\pr_R \colon V(\G) \to V_R(\G), \, \, a \mapsto \pr_R(a) \co a_r$ is a graph epimorphism from $\Gamma(\G)$ onto $\Gamma_R(\G)$. Consequently, if $\Gamma(\G)$ is connected, then so is $\Gamma_R(\G)$.
        \item\label{Prop4.5-ii} For any distinct vertices $u, v \in V_R(\G)$ and every $z_1, z_2 \in \Zz(\G)$, we always have
        \[\bigl(\{u, v\} \in E_R(\G)\bigr) \iff
        \bigl(\{u + z_1, v + z_2\} \in E(\G)\bigr).
        \]
		As a result, both $\Gamma_R(\G)$ and $\Gamma(\G)$ are necessarily incomplete; moreover, every edge in $\Gamma_R(\G)$ can be lifted to an edge in $\Gamma(\G)$ and the preimage under $\pr_R$ of any connected subgraph of $\Gamma_R(\G)$ is connected in $\Gamma(\G)$. 
	\end{enumerate}
\end{proposition}

\begin{proof}
We will prove each assertion separately. 
\begin{enumerate}[(i)] 
	\item {\bf Proof of~\eqref{Prop4.5-i}.} Since $\G = \G_R \oplus \Zz(\G)$ as a direct sum of vector spaces, for every $a \in V(\G) = \G \setminus \Zz(\G) \subseteq \G$, there exist unique $a_r \in \G_R$ and $a_z \in \Zz(\G)$ such that $a = a_r + a_z$. Since $a \notin \Zz(\G)$, we must have $a_r \neq 0$, hence $a_r \in (\G_{R})_0 = V_R(\G)$.
    \item[] Assume that $a, b \in V(\G); \, a_r, b_r \in V_R(\G); \, a_z, b_z \in \Zz(\G)$ with $a \ne b, \, a = a_r + a_z$ and $b = b_r + b_z$. By calculating directly and using the bilinearity of the Lie bracket, we obtain $[a, b] = [a_r + a_z, b_r + b_z] = [a_r, b_r]$.     
    It follows that
    \begin{align*}
    \bigl(\{a, b\} \in E(\G)\bigr) &\iff    
    (a \ne b) \wedge \bigl([a, b] = [a_r, b_r] = 0\bigr),\\
    &\iff (a_r = b_r) \vee \bigl(\{a_r, b_r\} \in E_R(\G)\bigr),
    \end{align*}
    which is exactly Relation~\eqref{For4.1New}. 
    \item[] Now, we consider the natural linear projection $\pr \colon \G \to \G_R$ associated to the decomposition $\G = \G_R \oplus \Zz(\G)$. Since $V(\G) = \G \setminus \Zz(\G) \subseteq \G$, we can consider the restriction $\pr_R \co \pr|_{V(\G)}$. Namely, thanks to the expansion in Item~\eqref{Prop4.5-i} just proven above, we have
    \[
		\pr_R = \pr|_{V(\G)} \colon V(\G) \to V_R(\G), \, \,  a \mapsto \pr_R(a) \co a_r \quad \text{where} \, \, a - a_r \in \Zz(\G).
	\]     
    Clearly $\pr_R$ is surjective because $\pr_R(v) = v$ for every $v \in V_R(\G) \subseteq V(\G)$. It remains to show that $\pr_R$ is a full homomorphism as the condition in Definition~\ref{Def1.2New}. Let $a \neq b$ be vertices of $\Gamma(\G)$. Since $a_r = \pr_R(a), \, b_r = \pr_R(b) \in V_R(\G)$, Relation~\eqref{For4.1New} confirms that the condition for a graph full homomorphism in the sense of Definition~\ref{Def1.2New} is satisfied. Thus, the surjective map $\pr_R\colon \Gamma(\G) \to \Gamma_R(\G)$ is indeed a graph epimorphism. Consequently, if $\Gamma(\G)$ is connected, then so is $\Gamma_R(\G)$. So, Assertion~\eqref{Prop4.5-i} holds.
    
	  \item {\bf Proof of~\eqref{Prop4.5-ii}.} Given any $u, v \in V_R(\G)$ ($u \ne v$) and any $z_1, z_2 \in \Zz(\G)$,  it is clear that $u + z_1 \ne v + z_2$, because if $u + z_1 = v + z_2$, then $u - v = z_2 - z_1 \in \G_R \cap \Zz(\G) = \{\boldsymbol{0}\}$, which contradicts the assumption $u \ne v$. Furthermore, $u + z_1 \in \pr_R^{-1}(u) \subseteq V(\G)$ and $v + z_2 \in \pr_R^{-1}(v) \subseteq V(\G)$. From the equality $[u, v] = [u+z_1, v+z_2]$, we immediately deduce    
      \[
      \bigl(\{u, v\} \in E_R(\G)\bigr) \iff
        \bigl(\{u + z_1, v + z_2\} \in E(\G)\bigr).
      \]
      Now, if $\Gamma_R(\G)$ is an LSG-complete, then $\G_R$ must be an abelian subalgebra of $\G$; in particular, $\G = \G_R \oplus_{\text{Lie}} \Zz(\G)$ is also abelian, which contradicts the assumption. In other words, $\Gamma_R(\G)$ and $\Gamma(\G)$ cannot be complete graphs.
      Specifically, any edge $e_r = \{u, v\}$ of $\Gamma_R(\G)$ can be lifted (in general, in many, even infinitely many, ways) to an edge $e = \{u+z_1, v+z_2\}$ of $\Gamma(\G)$. As a result, the preimage under $\pr_R$ of any connected subgraph of $\Gamma_R(\G)$ is connected in $\Gamma(\G)$. Hence, Assertion~\eqref{Prop4.5-ii} holds.      
	\end{enumerate}
The proof of Proposition~\ref{Prop4.5New} is complete.		
\end{proof}

\begin{theorem}[{\bf RCG as a useful version of CG}]\label{thm4}
 
Let $\G$ be a finite-dimensional non-abelian Lie algebra and $\Gamma(\G), \Gamma_R(\G)$ be the CG, RCG of $\G$, respectively. Let $\pr_R \colon V(\G) \to V_R(\G); \, a \mapsto \pr_R(a) = a_r$, be the graph epimorphism as in Item~\eqref{Prop4.5-i} of Proposition~\ref{Prop4.5New}. In that case, we have the following statements.
	\begin{enumerate}[(i)]
		\item\label{thm4-i} Let $\dist_\G \colon V(\G) \times V(\G) \to \N \cup \{+\infty\}$ and $\dist_R \colon V_R(\G) \times V_R(\G) \to \N \cup \{+\infty\}$ be the distance functions on $\Gamma(\G)$ and $\Gamma_R(\G)$, respectively. Then we have
        \begin{equation}\label{For4.2New}
        (a, b) \bigl(\in V(\G) \times V(\G)\bigr) \mapsto \dist_\G(a,b)=
            \begin{cases}
              0, & \text{if } \, \, a = b\\
	           1, & \text{if } \, \, a\ne b, \, a_r=b_r\\
	           \dist_R(a_r,b_r), & \text{if } \, \, a\ne b, \, a_r\neq b_r,
            \end{cases}        
        \end{equation}        
        Consequently, $\Gamma(\G)$ is connected if and only if $\Gamma_{R}(\G)$ is connected. 
        
		\item\label{thm4-ii} If $\Gamma(\G)$ and $\Gamma_{R}(\G)$ are both connected, their diameters are both infinite or both finite and equal. Let $n \co \dim\bigl(\Gamma_R(\G)\bigr)=\dim(\G) - \dim\bigl(\Zz(\G)\bigr)$. Then we always have 
        \begin{equation}\label{For4.3New}
            2 \leq h_{LSG}\bigl(\Gamma_R(\G)\bigr) \leq \min\Bigl{\{}n-1, \, \diam\bigl(\Gamma(\G)\bigr)\Bigr{\}}, \, n \geq 3
        \end{equation}
        
		\item\label{thm4-iii} If $\Gamma(\G)$ and $\Gamma_{R}(\G)$ are both disconnected, their sets of connected components have the same cardinality.
	\end{enumerate}
\end{theorem}

\begin{proof}[{\bf Proof of Theorem~\ref{thm4}}]

We will prove each statement of Theorem~\ref{thm4}.
\begin{itemize}
    \item {\bf Proof of~\eqref{thm4-i}.} 
    First, we will demonstrate Relation~\eqref{For4.2New} concerning the relationship between two distance functions $\dist_{\G}, \, \dist_R$. Let $a = a_r + a_z, \, b = b_r + b_z$ 
    be distinct vertices of $\Gamma(\G)$ with $a_r, b_r \in \G_R$, $a_z, b_z \in \Zz(\G)$. By Relation~\eqref{For4.1New} in Proposition~\ref{Prop4.5New}, we always have $[a, b] = [a_r, b_r]$. We need to compare the distance $\dist_{\G}(a, b)$ with $\dist_R(a_r, b_r)$. 
    \begin{itemize}
        \item[$\ast$] {\bf Case 1: $\boldsymbol{a_r = b_r}$:} $[a, b] = [a_r, b_r] = 0$. Therefore, $\{a, b\} \in E(\G)$. In particular, $\dist_R(a_r, b_r) = 0, \, \dist_{\G}(a, b) = 1$.
        \item[$\ast$] {\bf Case 2: $\boldsymbol{a_r \ne b_r}$:} Then, $\dist_{\G}(a, b) \geq 1$ and $\dist_R(a_r, b_r) \geq 1$. We will show that $\dist_{\G}(a, b) = \dist_R(a_r, b_r)$.
        \begin{itemize}
            \item[$\circ$] Prove $\dist_{\G}(a, b) \leq \dist_R(a_r, b_r)$: If $\dist_R(a_r, b_r) = +\infty$, the inequality is trivial. Suppose $\dist_R(a_r, b_r) = d < +\infty$ and 
            $a_r=v_0,v_1,\dots,v_d=b_r$ is a path of length $d$ connecting $a_r$ and $b_r$ in $\Gamma_R(\G)$. Then $a,v_1,\dots,v_{d-1},b$ is a path of the same length $d$ connecting $a$ and $b$ in $\Gamma(\G)$. Therefore,  $\dist_{\G}(a, b) \leq d = \dist_R(a_r, b_r)$. 
            
            \item[$\circ$] Prove $\dist_R(a_r, b_r) \leq \dist_{\G}(a, b)$: If $\dist_{\G}(a, b) = +\infty$, the inequality is trivial. Suppose $\dist_{\G}(a, b) = d < +\infty$ and 
            $a=u_0,u_1,\dots,u_d=b$ is a path of the length $d$ connecting $a$ and $b$ in $\Gamma(\G)$. Then, by Relation~\eqref{For4.1New} in Proposition~\ref{Prop4.5New}, any two consecutive terms of $a_r = \pr_R(a) = \pr_R(u_0), \pr_R(u_1), \cdots, \pr_R(u_d) = \pr_R(b) = b_r$ are equal or adjacent in $\Gamma_R(\G)$. Deleting consecutive repetitions, we obtain a path from $a_r$ to $b_r$ in $\Gamma_R(\G)$ of length at most $d = \dist_{\G}(a, b)$. Therefore, $\dist_R(a_r, b_r) \leq \dist_{\G}(a, b)$.
            \item[] In other words, $\dist_{\G}(a, b) = \dist_R(a_r, b_r)$. 
        \end{itemize} 
    \end{itemize} 
    Thus, Relation~\eqref{For4.2New} holds. Now, we will demonstrate that the connectedness of $\Gamma(\G)$ is equivalent to that of $\Gamma_R(\G)$.
    \begin{itemize}
        \item[$\ast$] If $\Gamma(\G)$ is connected, the graph epimorphism $\pr_R \colon V(\G) \to V_R(\G)$ implies that $\Gamma_R(\G)$ is also connected.
        
        \item[$\ast$] Conversely, if $\Gamma_R(\G)$ is connected, then the right-hand side of Relation~\eqref{For4.2New} is finite for all distinct $a,b\in V(\G)$. Thus, $\Gamma(\G)$ is connected.        
    \end{itemize}
    \item {\bf Proof of~\eqref{thm4-ii}}. Assume that both $\Gamma(\G)$ and $\Gamma_R(\G)$ are connected with diameters $\diam (\Gamma(\G))$ and $\diam (\Gamma_R(\G))$, respectively. We need to prove that their diameters are either both infinite or both finite and equal. Define
	\[
		\begin{array}{l}
			\dist(\G) \co \bigl{\{} d \in \N_0 \mid \exists \, a, b \in V(\G), a \neq b, \dist_\G(a, b) = d \bigr{\}}, \\
			\dist(R) \co \bigl{\{} d \in \N_0 \mid \exists \, u, v \in V_R(\G), u \neq v, \dist_R(u, v) = d \bigr{\}}.
		\end{array}
	\]	       
	We will prove that $\dist(\G)= \dist(R)$. First of all, by Relation~\eqref{For4.2New}, we have $\dist(R) \subseteq \dist(\G) \subseteq \{1\} \cup \dist(R)$. On the other hand, for any $u \in V_R(\G)$, every $\alpha \in \F_0$ such that $u \ne \alpha u$, we have $\dist_R(u, \alpha u) = 1 \in \dist(R)$. That means $\{1\} \cup \dist(R) \subseteq \dist(R)$. Hence $\dist(\G) \equiv \dist(R)$. By Definition~\ref{Def1.1New}(9), we get $\diam\bigl(\Gamma(\G)\bigr)=\diam\bigl(\Gamma_R(\G)\bigr)$ ($\in \N\cup\{+\infty\}$). 
    \item[] To prove Relation~\eqref{For4.3New}, put $n\co\dim\G_R=\dim(\G)-\dim\bigl(\Zz(\G)\bigr)$. By assumption, $\Gamma_R(\G)$ is an $n$-LSG which is connected and not complete. In particular, $\diam\bigl(\Gamma_R(\G)\bigr) \geq 2$. Since $\diam\bigl(\Gamma(\G)\bigr) = \diam\big(\Gamma_R(\G)\bigr) \geq 2$, $\Gamma(\G)$ must necessarily be incomplete.    
    Now, applying Theorem~\ref{thm1}~\eqref{thm1-iii} with $\diam\bigl(\Gamma_R(\G)\bigr) = \diam\bigl(\Gamma(\G)\bigr)$, we obtain $n \geq 3$ and Relation~\eqref{For4.3New}. 
    
    \item {\bf Proof of~\eqref{thm4-iii}}. Assume that both $\Gamma(\G)$ and $\Gamma_R(\G)$ are disconnected.	    
    Denote by $CC_\G(a)$ and $CC_R(v)$ the connected components of $a$ in $\Gamma(\G)$ and of $v$ in $\Gamma_R(\G)$, respectively. Let $a = a_r + a_z, \, b = b_r + b_z \in V(\G)$ with $a_r, b_r \in \G_R$, $a_z, \, b_z \in \Zz(\G)$. By Relation~\eqref{For4.2New}, $b$ is connected to $a$ in $\Gamma(\G)$ if and only if $b_r=a_r$ or $b_r$ is connected to $a_r$ in $\Gamma_R(\G)$. That means, for any $a \in V(\G)$
    \begin{equation}\label{For4.4New}
        CC_{\G}(a_r) = CC_\G(a)=\pr_R^{-1}\bigl(CC_R(a_r)\bigr); \, \pr_R\bigl(CC_{\G}(a)\bigr) = CC_R(a_r)
    \end{equation}
    Denote by 
	\[
		CC_{\Gamma(\G)} \co \bigl{\{}CC_{\G}(a) \mid a \in V(\G)\bigl{\}}; \, \, CC_{\Gamma_R(\G)} \co \bigl{\{}CC_R(v) \mid v \in V_R(\G)\bigl{\}}
	\]
    the sets of connected components of $\Gamma(\G)$ and $\Gamma_R(\G)$, respectively. We will prove
	\[
		\text{card} \bigl(CC_{\Gamma(\G)}\bigr) = \text{card} \bigl(CC_{\Gamma_R(\G)}\bigr)
	\]
    Relation~\eqref{For4.4New} gives us the two mappings below
    \begin{alignat*}{2}
    \Phi & \colon CC_{\Gamma(\G)} \to CC_{\Gamma_R(\G)}; \quad CC_{\G}(a) \mapsto \Phi\bigl(CC_{\G}(a)\bigr) \co CC_R(a_r);\\
    \Psi &\colon CC_{\Gamma_R(\G)} \to CC_{\Gamma(\G)}; \quad 
     CC_R(a_r) \mapsto \Psi\bigl(CC_R(a_r)\bigr) \co CC_{\G}(a_r).
    \end{alignat*}
    Furthermore, thanks to the surjectivity of $\pr_R$ and Relation~\eqref{For4.4New}, we can verify that $\Phi$ and $\Psi$ are inverses of each other, and are therefore bijections. Thus, $\text{card} \bigl(CC_{\Gamma(\G)}\bigr) = \text{card} \bigl(CC_{\Gamma_R(\G)}\bigr)$. 
    \end{itemize} Theorem~\ref{thm4} is fully proven. 
\end{proof}

\begin{theorem}[{{\bf Decomposability of the RCG}}]\label{thm5}
	Let $\G=\G_1\oplus_{\text{Lie}}\G_2$ be a direct sum of two non-abelian Lie algebras. Then $\Gamma_R(\G_1\oplus_{\text{Lie}}\G_2)\cong_{LSG}\Gamma_R(\G_1)\oplus_{LSG}\Gamma_R(\G_2)$.
\end{theorem} 

\begin{proof}[{\bf Proof of Theorem~\ref{thm5}}]
We have $\Zz(\G)=\Zz(\G_1)\oplus\Zz(\G_2)$. By Theorem~\ref{thm3} we may take $\G_R=\G_{1R}\oplus\G_{2R}$, where $\G_{iR}$ is a complement of $\Zz(\G_i)$ in $\G_i$. For $x=(x_1,x_2)$, $y=(y_1,y_2)\in (\G_R)_0$, we have
\[
\{x,y\}\in E_R(\G) \iff \bigl(x_1\sim_{LSG}y_1 \, \, \text{ in }\Gamma_R(\G_1)\bigr) \wedge \bigl( x_2\sim_{LSG}y_2\, \, \text{in }\Gamma_R(\G_2)\bigr).
\]
By Definition~\ref{Def2.5New} and Proposition~\ref{Prop2.6New}, the identity map $\id \colon \G_R \overset{\cong}{\to} \G_{1R} \oplus \G_{2R}$ induces an LSG-isomorphism
\[LSG(\id)\colon \Gamma_R(\G_1 \oplus \G_2) \overset{\cong}{\to}_{LSG} \Gamma_R(\G_1) \oplus_{LSG} \Gamma_R(\G_2).\]
Theorem~\ref{thm5} is fully proven.
\end{proof}

By Theorem~\ref{thm5} it suffices to consider RCGs of Lie algebras that do not decompose into two non-abelian summands. In the rest of the paper we describe explicitly the RCGs (and hence, by Theorem~\ref{thm4}, the key data for the CGs) of all AALAs.

 
\subsection{Reduced commuting graphs of AC Lie algebras}\label{Subsec4.3New}
AC Lie algebras were introduced by Sriwongsa~\cite{Sriw2025} in 2025.

\begin{definition}[{\cite{Sriw2025}}]\label{Def4.6New}
    A Lie algebra $\G$ is called an \emph{AC Lie algebra} (ACLA, for short) if  Relation~\eqref{For3.1New} in Definition~\ref{Def3.4New} is satisfied for every triple $x, y, z \in \G\setminus \Zz(\G)$, i.e.
    \[\bigl[([x, y]= 0) \wedge ([y, z] = 0)\bigr] \Longrightarrow ([x, z] = 0); \, \forall x, y, z \in \G\setminus \Zz(\G).\]
\end{definition}

\begin{remark}\label{Re4.7ACNew}
    Note that each CTLA $\G$ is always an ACLA, but converse is not true (see Example 3.11 in~\cite{Sriw2025}).
\end{remark}

The following proposition, which characterizes the class of AC Lie algebras, extends Theorem~\ref{thm2} to Lie algebras with an arbitrary center. Its proof follows the same lines as the proof of Theorem~\ref{thm2}, combined with~\cite[Lemma~3.12]{Sriw2025}, and is therefore omitted.

\begin{proposition}[{\bf LSG-height and RCGs of ACLAs}]\label{Prop4.8New}
   Let $\G$ be an $(n+m)$-dimensional non-abelian Lie algebra with an $m$-dimensional center, where $m,n\in\N$ and $n\geq 1$. Fix a complementary subspace $\G_R$ to the center in $\G$. Of course, $\dim(\G_R) = n$. Let $\Gamma_R(\G)$ be $n$-RCG of $\G$. The following are equivalent:
\begin{enumerate}[(i)]
    \item\label{Prop4.8-i} $h_{LSG}\bigl(\Gamma_R(\G)\bigr) \leq 1$.
    
    \item\label{Prop4.8-ii} $C_{\G}(s)$ is an abelian subalgebra of $\G$ for every $s \in (\G_{R})_0$.
    
    \item\label{Prop4.8-iii} $\G$ is an ACLA.
\end{enumerate}
If these conditions hold, then, for every $s \in (\G_R)_0$, the $C_{\G}(s)$ are exactly the maximal abelian subalgebras of $\G$, $\s_{LSG}(s) = C_{\G}(s) \cap \G_R$ and $CC_R(s) = \bigl(C_{\G}(s) \cap \G_R\bigr)_0$, where $CC_R(s)$ is the connected component of $s$ in $\Gamma_R(\G)$. Furthermore, $\Gamma_R(\G)$ is always disconnected and the LSG-completes $\boldsymbol{K}_{\bigl(C_{\G}(s) \cap \G_R\bigr)_0}$ are precisely the LSsubG-connected components of $\Gamma_R(\G)$, and 
\[
\bigl[h_{LSG}\bigl(\Gamma_R(\G)\bigr) = 0\bigr] \Leftrightarrow \bigl{\{}[C_{\G}(s) = \Zz(\G) \oplus \F s; \, \forall s \in (\G_{R})_0] \wedge \bigl[\Gamma_R(\G)\cong_{LSG}\bigsqcup_{\widetilde{s} \in \mathbb{P}^{n-1}_{\F}} \boldsymbol{K}_{\widetilde{s}}\bigr]\bigr{\}}
\]
\hfill $\square$
\end{proposition}

\subsection{Reduced commuting graphs of almost abelian Lie algebras}\label{Subsec4.4New}

Almost abelian Lie algebras were introduced by Kolman~\cite{Kol1965} in the classification of upper semimodular Lie algebras, later extended by Gein~\cite{Gei76}. We adopt the broad definition below.

\begin{definition}[{\cite[Definition 1]{Ave22}}]\label{Def4.7New}
    A non-abelian Lie algebra $\G$ is an \emph{almost abelian Lie algebra} (AALA) if it contains an abelian subalgebra of codimension $1$. For an $(n+1)$-dimensional AALA, we often refer to it as an $(n+1)$-AALA.
\end{definition}

\begin{remark}\label{Re4.7New} We have a few observations below.
    \begin{enumerate}[(1)]
        \item\label{Re4.7-1} By~\cite[Proposition~3.1]{BC12}, an abelian subalgebra of codimension $1$ yields an abelian ideal of codimension $1$. So, every AALA is $2$-step solvable.
        
        \item\label{Re4.7-2} An $(n+1)$-AALA $\G$ always has the form $\G = \F^n \oplus_D \F e_0$ with $\boldsymbol{0}\neq D \in \text{Der}(\F^n) \equiv \mathfrak{gl}_n(\F)$ and $\Zz(\G) = \Ker(D)$.
        
        \item\label{Re4.7-3} By Item~\eqref{Re4.7-2}, $\Zz(\G)=\Ker(D)$, so an $(n+1)$-AALA with trivial center has $D\in GL_n(\F)$; this applies in particular to cocomplete AALAs by~\cite[Theorem 3]{Vu2025}.
    \end{enumerate}
\end{remark}

\begin{notation}\label{No4.8New}
Let $\F^{n+1}\equiv\F^n\oplus\F e_0$ with $\F^n\equiv\F^n\times\{0\}$. For $s\in\F^{n+1}\setminus\F^n$, $\boldsymbol{K}_{\widetilde{s}}$ is the complete $1$-LSG on $\F_0 s$. An $(n+m+1)$-AALA with $m$-dimensional center is written
\begin{equation}\label{For4.5New}
   \G = \F^{n+m} \oplus _D \F e_0; \, \, e_0 \notin \F^{n+m},\, \, D \in \mathfrak{gl}_{n+m}(\F),\ [e_0,v]=D(v),\,\, \Zz(\G) = \Ker (D)
\end{equation}
\end{notation}

We focus on AALAs for three reasons: they are the non-abelian Lie algebras closest to abelian ones and play a role in the classification of solvable Lie algebras~\cite{Ave22, BC12, Vu2022, Vu2025}; the class is large, with centers of arbitrary dimension, so the redundancy removed by the RCG appears in full generality; and their RCGs have LSG-height $0$ or $1$, exactly the range in which Theorem~\ref{thm1} determines all components. Theorem~\ref{thm6} answers the main question in both directions: it computes the LSG-invariants of the RCG of an AALA, and shows conversely that these invariants force a Lie algebra $\G'$ to be an AALA, with $\dim\Zz(\G')=\dim\G'-\dim\Gamma_R(\G')$.

\begin{theorem}[{\bf RCGs of $(n+m+1)$-AALAs and their characterization}]\label{thm6}
Let $n \in \N_0$ and $m \in \N$.
\begin{enumerate}[(i)]
    \item\label{thm6-i} Let $\G$ be an $(n+m+1)$-AALA with $m$-dimensional center, written as in Relation~\eqref{For4.5New}, with $\Zz(\G)=\Ker(D) = \{\boldsymbol{0}\} \times \F^m$ and complement $\G_R \co \F^n \oplus \F e_0$. Then for every $x \in V_R(\G)$,
    \[
        \s_{LSG}(x) = N_{LSG}(x) =
        \begin{cases}
            \F^n , & x \in \F^n_0 ,\\
            \F x , & x \notin \F^n .
        \end{cases}
    \]
    In particular $h_{LSG}\bigl(\Gamma_R(\G)\bigr) = \min\{1,n-1\}$, $\dim_{LSG}\bigl(\Gamma_R(\G)\bigr) = n = \dim \Gamma_R(\G) - 1$, and
    \begin{equation}\label{For4.6New}
        \Gamma_R(\G) = \boldsymbol{K}_{\F^n_0} \ \sqcup \Bigl(\bigsqcup_{\widetilde{s} \in (\F^{n+1}_0 \setminus \F^n_0)/\sim \, =\mathbb{P}^{n}_{\F} \setminus \mathbb{P}^{n-1}_{\F}} \boldsymbol{K}_{\widetilde{s}} \Bigr) .
    \end{equation}
    So $\Gamma_R(\G)$ is disconnected, its components are the LSG-complete $\boldsymbol{K}_{\F^n_0}$ and the LSG-completenes $\boldsymbol{K}_{\widetilde{s}}$, $s \notin \F^n$.
    
    \item\label{thm6-ii} Conversely, let $\G'$ be an $(n+m+1)$-dimensional non-abelian Lie algebra with
    \[
        h_{LSG}\bigl(\Gamma_R(\G')\bigr) \leq 1, \quad  \dim_{LSG}\bigl(\Gamma_R(\G')\bigr) = \dim \Gamma_R(\G') - 1 = n .
    \]
    Then $\Gamma_R(\G')$ is of the form~\eqref{For4.6New} and $\G'$ is an $(n+m+1)$-AALA with $m$-dimensional center. Moreover, $h_{LSG}\bigl(\Gamma_R(\G')\bigr) = 0$ if and only if $n = 1$.
\end{enumerate}
\end{theorem}

\begin{proof}[{\bf Proof of Theorem~\ref{thm6}}]
We will prove each assertion.
\begin{itemize} 
    \item {\bf Proof of~\eqref{thm6-i}}: For $x = u + \alpha e_0$, $y = v + \beta e_0$ in $\G_R$ ($u, v \in \F^n$), $[x, y] = D(\alpha v - \beta u)$, and $D$ is injective on $\F^n$ since $\F^n \cap \Ker(D) = \{\boldsymbol{0}\}$. Hence $[x, y] = 0$ iff $\alpha v = \beta u$. If $x \in \F^n_0$ ($\alpha = 0$), this means $\beta = 0$, i.e.\ $y\in\F^n$; if $x \notin \F^n$ ($\alpha \neq 0$), it means $y = \alpha^{-1}\beta x \in \F x$. This gives $N_{LSG}(x)$. All non-zero vectors of $N_{LSG}(x)$ have the same LSG-neighborhood, so $N_{LSG}(x)$ is LSG-closed and $\s_{LSG}(x)=N_{LSG}(x)$. Moreover $h_{LSG}(x)=1$ iff $\F x\subsetneq N_{LSG}(x)$ if and only if ($x\in\F^n_0$ and $n\ge2$), and $h
    _{LSG}(x)=0$ otherwise; this gives the height and the closed dimension. Finally,  Relation~\eqref{For4.6New} is derived from Relation~\eqref{For2.5New} applied to the $(n+1)$-LSG $\Gamma_R(\G)$ when $n \geq 2$ or Relation~\eqref{For2.4New} of Theorem~\ref{thm1}(ii)(a) when $n = 1$.     
    \item {\bf Proof of~\eqref{thm6-ii}}: 
    Put $\Gamma \co \Gamma_R(\G')$, constructed on a complement $\G'_R$ of $\Zz(\G')$. Since $\dim\G'_R=n+1$, $\dim \Zz(\G') = m$. We find an $n$-dimensional subspace $W \leq \G'_R$ with $W_0$ a clique and $\Gamma$ of the form Relation~\eqref{For4.6New}.
    \begin{itemize}
    \item If $h_{LSG}(\Gamma) = 0$, then $N_{LSG}(s) = \F s$ for all $s$, so $n = \dim_{LSG}(\Gamma) = 1$; by Relation~\eqref{For2.4New}, $\Gamma = \bigsqcup_{\widetilde{s} \in \mathbb{P}^{1}_{\F}} \boldsymbol{K}_{\widetilde{s}}$, which is Relation~\eqref{For4.6New} for $n = 1$. Then we take $W=\F s_0$ for any $s_0\in(\G'_R)_0$.
    \item If $h_{LSG}(\Gamma) = 1$, then $n = \dim_{LSG}(\Gamma) \geq 2$ by Remark~\ref{Re2.15New}, and by Theorem~\ref{thm1}(ii)(b) there is $s_0$ with $W \co N_{LSG}(s_0)$ of dimension $n = \dim \Gamma - 1$, $\boldsymbol{K}_{W_0}$ a connected component, and $\Gamma$ given by Relation~\eqref{For2.5New}.     
    \end{itemize}
    \item[] In both cases, one has $[W, W] = \{\boldsymbol{0}\}$. Since $\Zz(\G')$ is central, $A \co W \oplus \Zz(\G')$ is an abelian subalgebra of dimension $n + m = \dim(\G') - 1$ and $\G'$ is an AALA of dimension $(n+m+1)$ with $m$-dimensional center. The last claim follows from the proof of~\eqref{thm6-i}, since both conditions can be determined from Relation~\eqref{For4.6New}.
\end{itemize} Theorem~\ref{thm6} is fully proven.
\end{proof}

\begin{remark}\label{Re4.9New}
For $m=0$, Theorem~\ref{thm6} provides the LAGs for all $(n + 1)$-AALAs with a trivial center, including the class of cocomplete AALAs
(see Remark~\ref{Re4.7New}(3)); their LAGs have the form~\eqref{For4.6New}, which generalizes Example~\ref{Ex3.3New}\eqref{Ex3.3-1}. In particular, $(n+1)$-AALAs with trivial center are recognized by their LAGs.
\end{remark}

\begin{corollary}\label{Co5.5New}
Any two AALAs of the same dimension whose centers have the same dimension have LSG-isomorphic RCGs. \hfill $\square$
\end{corollary}

\begin{remark}\label{Re5.6New}
The classes $\text{Lie}(n+i, n\text{C})$ ($i=1,2$) classified in~\cite[Theorem 1]{Vu2022} include the subclass $\text{Lie}(n+1, n\text{C})$ of $(n+1)$-AALAs, whose RCGs are thus described by Theorem~\ref{thm6} and Corollary~\ref{Co5.5New}. The next step is to describe the RCGs of $\text{Lie}(n+2, n\text{C})$, followed by $\text{Lie}(n+2, n)$.
\end{remark}

\section{Conclusion, open problems and declaration of AI use}\label{Sec5}

\subsection{Summary of the main results}\label{Subsec5.1New}
This paper continues the line of research 
in~\cite{HHVT2024}, guided by the recognition problem stated in the Introduction. We introduced LSGs and, as special cases, LAGs and RCGs. The RCG is a faithful replacement for the CG: it preserves the connectivity, the diameter and the number of connected components, removes the redundancy created by the center, and retains the linear structure. We characterized LSGs of LSG-height at most 1 and showed that an LAG has LSG-height at most 1 if and only if the Lie algebra is CT. We also provided a complete description of the RCGs of all AALAs. The class AALA is recognized by the RCGs; moreover, the RCG $\Gamma_R(\G)$ of any AALA $\G$ determines only $\dim(\G) - \dim \bigl(\Zz(\G)\bigr)$, therefore, once $\dim\bigl(\Zz(\G)\bigr)$ is known, $\dim(\G)$ is known. 

\subsection{Open problems}\label{Subsec5.2New} 
On the basis of this study of LSGs, LAGs and RCGs, we propose three open problems for further research.
 
First of all, recall that in 2006 Abdollahi, Akbari and Maimani~\cite{AAM2006} conjectured that every finite non-abelian simple group is determined by its non-commuting graph; this was proved by Solomon and Woldar~\cite{Solomon2013} in 2013. Inspired by this result, we consider an analogous question for Lie algebras and their RCGs. Theorem~\ref{thm6} shows that the RCG recognizes the whole class of AALAs, but it is not truly ``rigid'': two Lie algebras belonging to the AALA class may be non-isomorphic even if their RCGs are isomorphic. Example~\ref{Ex4.3New} exhibits two non-isomorphic Lie algebras with the same RCG. Therefore, we pose the first open problem to be addressed: characterizing and classifying the Lie algebras with trivial center whose LAGs possess a ``rigid'' property. Specifically, we present the first open problem as follows.

\begin{problem}\label{OP1New} 
Describe and classify (up to Lie algebra isomorphism) the finite-dimensional Lie algebras $\G$ with trivial center such that, for every finite-dimensional Lie algebra $\G'$ with trivial center, the following property holds: 
\begin{equation}\label{For5.1New}
    \bigl(\Gamma(\G) \cong_{LSG} \Gamma(\G')\bigr) \Longrightarrow \bigl(\G \cong_{\text{Lie}} \G'\bigr).
\end{equation}
\end{problem}

\begin{remark}\label{Re5.1New}
Regarding Open problem~\ref{OP1New}, with the assistance of Claude Opus 5.5 (Anthropic), we conjecture that, over an algebraically closed field $\F$ of characteristic zero, any two simple Lie algebras with LSG-isomorphic LAGs are isomorphic. Algebraic closure cannot be dropped: over $\R$, all non-zero elements of $\mathfrak{so}(3)$ and of $\mathfrak{sl}_2(\R)$ have $1$-dimensional centralizers, so $\Gamma(\mathfrak{so}(3)) \cong_{LSG} \Gamma(\mathfrak{sl}_2(\R))$, although $\mathfrak{so}(3)$ and $\mathfrak{sl}_2(\R)$ are not Lie isomorphic. However, identifying the class of solvable Lie algebras that satisfy Relation~\eqref{For5.1New} will likely prove to be a significant challenge.
\end{remark}

We now return to the Main question posed in the Introduction, reformulated for RCGs: which classes of Lie algebras are determined by the LSG-invariants of their RCGs? Theorems~\ref{thm2} and~\ref{thm6} describe two classes of Lie algebras whose RCGs have LSG-height at most $1$, namely CT Lie algebras and AALAs. On the other hand, by Theorem~\ref{thm1}(ii)(a),(b), for such $\G$ the RCG $\Gamma_R(\G)$ is a disjoint union of LSG-completes $\boldsymbol{K}_{W_0}$, where the subspaces $W = N_{LSG}(s)$ form a vector space partition of $\G_R$. Two cases are already settled in the paper. First, if $\Zz(\G) = \{0\}$, these are exactly the CTLAs (Theorem~\ref{thm2}), classified by Klep and Moravec~\cite{Klep2010} over algebraically closed fields of characteristic zero and by Gorbatsevich~\cite{Gorb2017} over $\R$. Second, if $\G$ is almost abelian, Condition~\eqref{For5.2New} always holds, whatever the dimension of the center (Theorem~\ref{thm6}). Furthermore, this is also linked to the class of AC Lie algebras recently introduced in \cite{Sriw2025} by Sriwongsa (2025). Note that the class of AC Lie algebras (ACLAs) is broader than the class of CTLAs and has not yet been classified. Furthermore, similar to Theorem~\ref{thm2}, Proposition~\ref{Prop4.8New} asserts that, for any ACLA $\G$, its RCG always has the height of at most 1. This leads us to the second open problem, presented below.

\begin{problem}[{{\bf Lie algebras whose RCG has LSG-height at most $1$}}]\label{OP2New}
Classify the finite-dimensional non-abelian Lie algebras $\G$ such that
\begin{equation}\label{For5.2New}
    h_{LSG}\bigl(\Gamma_R(\G)\bigr) \leq 1.
\end{equation}
Specifically, classify all non-abelian ACLAs. 
\end{problem}

On the other hand, the inequality $2 \leq h_{LSG}(\Gamma) \leq \min\bigl{\{}\dim(\Gamma) - 1, \, \diam(\Gamma)\bigr{\}}$ in Theorem~\ref{thm1}\eqref{thm1-iii}, together with  Relation~\eqref{For4.3New} in Theorem~\ref{thm4}, reveals a close connection between the height of an LSG (or RCG) and its diameter.
This relationship always holds true for LSGs (or RCGs) which are connected and not complete. Naturally, one might consider examining connected but not complete LSGs (or RCGs) for which the height is equal to the diameter. This is the subject of the third open problem presented below.

\begin{problem}[{{\bf Lie algebras whose RCG has LSG-height equal to its diameter}}]\label{OP3New}
Let $3 \leq n \in \N$.
\begin{enumerate}[(1)]
    \item\label{OP3-1} Characterize the connected $n$-LSGs $\Gamma$ constructed on $\F^n$ such that $h_{LSG}(\Gamma) = \diam(\Gamma) \geq 2$.
    
    \item\label{OP3-2} Classify the finite-dimensional Lie algebras $\G$ over $\F$ with $n = \dim(\G) - \dim\Zz(\G)$ such that the RCG and CG of $\G$ satisfy the property
    \[
    h_{LSG}\bigl(\Gamma_R(\G)\bigr) = \diam\bigl(\Gamma(\G)\bigr) \geq 2.
    \]    
\end{enumerate}
\end{problem}

\begin{remark}\label{Re5.2New} Let us return to 
  the $3$-LSG $\Gamma$ in Example~\ref{Ex2.17New} with $\F = \R$. A direct calculation gives    
  $h_{LSG}(\Gamma) = \diam(\Gamma) = 2$. Thus, $\Gamma$ is an example of the kind of LSG considered in Open problem~\ref{OP3New}\eqref{OP3-1}.  
\end{remark}

\subsection{Declaration of generative AI use}\label{Subsec5.3New}
\addcontentsline{toc}{section}{Declaration of generative AI use}

All definitions and results presented in this paper were conceived and developed by the authors. The initial proofs were also constructed by the authors and subsequently refined with the assistance of AI. Specifically, during the preparation of this manuscript, the first author used Claude Opus 5.5 (Anthropic), while the other authors used Claude Fable 5.1 (Anthropic) to assist with the following tasks:
\begin{enumerate}[(1)]
    \item English-language editing, LaTeX typesetting, and checking proofs and routine computations.
    \item Streamlining the presentation of the Abstract and the Introduction. Refining the statements and proofs of Theorem~\ref{thm1}\eqref{thm1-i}, Theorem~\ref{thm4}, Theorem~\ref{thm6}\eqref{thm6-ii}, and the presentation of Example~\ref{Ex2.18New}.
    \item Improving the clarity and correctness of the exposition.    
\end{enumerate}
It was Claude Opus 5.5 (from Anthropic) that drew the authors' attention to the concept of ``AC Lie algebras'' in~\cite{Sriw2025}. Consequently, the authors defined a research task concerning the class of ACLAs within the context of Open problem~\ref{OP2New}.
Every passage produced with the assistance of these tools was checked, revised and verified by the authors. We take full responsibility for the content of this article.

\section*{Disclosure statement}
\addcontentsline{toc}{section}{Disclosure statement}
No potential conflict of interest was reported by the author(s).


\section*{Funding}
\addcontentsline{toc}{section}{Funding}
This research is funded by Vietnam National University Ho Chi Minh City (VNU-HCM) under grant number B2026-34-06. 


\end{document}